\documentclass[12pt]{article}
\usepackage{graphicx}
\usepackage{amsmath,amsthm,amssymb,amsfonts,euscript,enumerate}
\usepackage{amsthm}
\usepackage{color,wrapfig}
\usepackage{pxfonts}
\usepackage{verbatim}
\newtheorem{thm}{Theorem}[section]
\newtheorem{cor}[thm]{Corollary}
\newtheorem{lem}[thm]{Lemma}

\newtheorem{defn}[thm]{Definition}
\newtheorem{rmk}[thm]{Remark}

\newcommand{\Rr}{{\mathbb{R}}}
\newcommand{\Z}{{\mathbb{Z}}}

\newcommand{\1}{\partial}
\newcommand{\2}{\overline}
\newcommand{\3}{\varepsilon}

\begin{document}
\title{Asymptotic behaviour of the finite  blow-up points solutions of the fast diffusion equation}
\author{Shu-Yu Hsu\\
Department of Mathematics\\
National Chung Cheng University\\
168 University Road, Min-Hsiung\\
Chia-Yi 621, Taiwan, R.O.C.\\
e-mail: shuyu.sy@gmail.com}
\date{Aug 5, 2023}
\smallbreak \maketitle
\begin{abstract}
Let $n\ge 3$, $0<m<\frac{n-2}{n}$, $i_0\in\Z^+$, $\Omega\subset\mathbb{R}^n$ be a smooth bounded domain, $a_1,a_2,\dots,a_{i_0}\in\Omega$, $\widehat{\Omega}=\Omega\setminus\{a_1,a_2,\dots,a_{i_0}\}$,  $0\le f\in L^{\infty}(\partial\Omega)$  and $0\le u_0\in L_{loc}^p(\widehat{\Omega})$ for some constant $p>\frac{n(1-m)}{2}$ which satisfies $\lambda_i|x-a_i|^{-\gamma_i}\le u_0(x)\le \lambda_i'|x-a_i|^{-\gamma_i'}\,\,\forall 0<|x-a_i|<\delta$, $i=1,\dots, i_0$ where $\delta>0$, $\lambda_i'\ge\lambda_i>0$ and $\frac{2}{1-m}<\gamma_i\le\gamma_i'<\frac{n-2}{m}$ $\forall i=1,2,\dots, i_0$  are constants. We will prove the asymptotic behaviour of the finite blow-up points solution $u$ of $u_t=\Delta u^m$ in $\widehat{\Omega}\times (0,\infty)$, $u(a_i,t)=\infty\,\,\forall i=1,\dots,i_0, t>0$, $u(x,0)=u_0(x)$ in $\widehat{\Omega}$ and $u=f$ on $\partial\Omega\times (0,\infty)$, as $t\to\infty$. We will construct  finite blow-up points  solution in bounded cylindrical domain with appropriate  lateral boundary value such that the finite blow-up points  solution oscillates between two given harmonic functions as $t\to\infty$.  We will also prove the existence of the minimal  solution of $u_t=\Delta u^m$ in $\widehat{\Omega}\times (0,\infty)$, $u(x,0)=u_0(x)$ in $\widehat{\Omega}$, $u(a_i,t)=\infty\quad\forall t>0, i=1,2\dots,i_0$ and $u=\infty$ on $\partial{\Omega}\times (0,\infty)$.
\end{abstract}

\vskip 0.2truein



\vskip 0.2truein
\setcounter{section}{0}

\section{Introduction}
\setcounter{equation}{0}
\setcounter{thm}{0}

Recently there is a lot of study on the properties of the fast diffusion equation,
\begin{equation}\label{fde}
u_t=\Delta u^m
\end{equation}   
with $0<m<1$ by P.~Daskalopoulos,  M.~Fila, S.Y.~Hsu, K.M.~Hui, T.~Jin, S.~Kim, Y.C.~Kwong, P.~Mackov\'a, M. del Pino, M.~S\'aez,  N.~Sesum, J.~Takahashi, J.L. Vazquez, H.~Yamamoto, E.~Yanagida, M.~Winkler and J.~Xiong, etc. \cite{DS1}, \cite{DS2}, \cite{FMTY}, \cite{H1}, \cite{H2}, \cite{HK1}, \cite{HK2}, \cite{Hs}, \cite{JX}, \cite{K}, \cite{PS}, \cite{TY}, \cite{VW1}, \cite{VW2}. The equation \eqref{fde} arises in many physical models and in geometry. When $m>1$, \eqref{fde} is called the porous medium equation which arises in the modelling of gases passing through porous media and oil passing through sand  \cite{Ar}. \eqref{fde} also arises as the diffusive limit for the generalized Carleman kinetic equation \cite{CL}, \cite{GS}, and as the large time asymptotic limit of the solution of the free boundary compressible Euler equation with damping \cite{LZ}. When $m=1$, it is the heat equation. When the dimension $n\ge 3$ and $m=\frac{n-2}{n+2}$, \eqref{fde} arises in the study of the Yamabe flow \cite{DS1}, \cite{DS2}, \cite{PS}. 

Various fundamental results in $\Rr^n$ for the equation \eqref{fde} are obtained recently by A.~Friedman and S.~Kamin \cite{FrK}, and M.~Bonfonte,   J.~Dolbeault, G.~Grillo, M. del Pino and J.L.~Vazquez \cite{BBDGV}, \cite{BGV1}, \cite{BV1}, \cite{BV2},  \cite{CaV}, \cite{PD}, \cite{V1}. Results for the equation \eqref{fde} in bounded domains are also obtained by D.G.~Aronson and L.A.~Peletier \cite{ArP}, B.E.J.~Dahlberg and  C.E.~Kenig \cite{DaK}, E.~Dibenedetto and Y.C.~Kwong \cite{DiK}, E.~Feireisl and F.~Simondon \cite{FeS}, M.~Bonfonte,  G.~Grillo and J.L.~Vazquez \cite{BGV2}, \cite{V2}. We refer the readers to the books by P.~Daskalopoulos and C.E.~Kenig   \cite{DK} and J.L. Vazquez \cite{V3}, \cite{V4}, for some recent results on the equation \eqref{fde}.

Let $\Omega\subset\mathbb{R}^n$ be a smooth bounded domain, $i_0\in\Z^+$, $a_1,a_2,\dots,a_{i_0}\in\Omega$ and 
\begin{equation}\label{omega-hat}
\widehat{\Omega}=\Omega\setminus\{a_1,a_2,\dots,a_{i_0}\}.
\end{equation}
Let $\delta_0=\frac{1}{3}\min_{1\leq i,j\leq i_0}\left(\textbf{dist}(a_i,\partial\Omega),\left|a_i-a_j\right|\right)$. For any $\delta>0$, let 
\begin{equation}\label{d-delta}
\Omega_{\delta}=\Omega\setminus\left(\cup_{i=1}^{i_0}B_{\delta}(a_i)\right)\quad\mbox{ and }\quad D_{\delta}=\{x\in\Omega_{\delta}:\mbox{dist}\, (x,\partial\Omega)>\delta\}.
\end{equation}
 For any $j\in\Z^+$, let 
\begin{equation}\label{ej}
E_j=\{x\in\widehat{\Omega}:\mbox{dist }(x,\partial\Omega)>1/j\} 
\end{equation}
and $j_0\in\Z^+$ be such that $j_0>1/\delta_0$. 

When $n\ge 3$ and $\frac{n-2}{n}<m<1$, existence of positive smooth solution of the Cauchy problem,
\begin{equation}\label{Cauchy-problem}
\left\{\begin{aligned}
u_t=&\Delta u^m\quad\mbox{ in }\mathbb{R}^n\times (0,\infty)\\
u(x,0)=&u_0(x)\quad\mbox{in }\mathbb{R}^n
\end{aligned}\right.
\end{equation}
for any $0\le u_0\in L^1_{loc}(\mathbb{R}^n), u_0\ne 0$, was proved by M.A.~Herrero and M.~Pierre in \cite{HP}. This implies that the diffusion for the solution of \eqref{Cauchy-problem} when $m<1$ must be very fast so that for any $t>0$ the solution $u(x,t)$ is positive everywhere even though the initial value $u_0$ may only have compact support in $\Rr^n$.

When $m>1$, the Barenblatt solution
\begin{equation*}
B(x,t)=t^{-\alpha}\left(C_1-\frac{\alpha(m-1)|x|^2}{2mnt^{\alpha/n}}\right)_+^{\frac{1}{m-1}},\quad x\in\mathbb{R}^n, t>0
\end{equation*}
where $\alpha=\frac{n}{n(m-1)+2}$ and $C_1>0$ is any constant is a solution of \eqref{fde} in $\Rr^n\times (0,\infty)$. 
Note that the Barenblatt solution $B(x,t)$ has compact support for any $t>0$ and the diffusion for the solution of \eqref{fde} when $m>1$ is slow. Hence the behaviour of the solution of \eqref{fde} is very different for $m<1$ and $m>1$. 

 Uniqueness of solutions of \eqref{Cauchy-problem} for the case $0<m<1$ and $n\ge 1$ was also proved in \cite{HP}. This result was later extended by G.~Grillo, M.~Muratori and F.~Punzo \cite{GMP}, to the uniqueness of the strong solution of the equation
\begin{equation*}
\left\{\begin{aligned}
u_t=&\Delta u^m\quad\,\mbox{ in }M\times (0,\infty)\\
u(x,0)=&u_0(x)\quad\mbox{ in }M
\end{aligned}\right.
\end{equation*}
on a Riemannian manifold $M$ whose Ricci curvatures satisfies some lower bound condition. Note that in the uniqueness theorems of  \cite{HP} and \cite{GMP} the solutions considered are strong 
solutions. That is the solution $u$ satisfies 
\begin{equation*}
u_t\in L^1_{loc}(\Rr^n\times (0,\infty))\quad\mbox{ in \cite{HP}}
\end{equation*}
and 
\begin{equation*}
u_t\in L^1_{loc}(M\times (0,\infty)) \quad\mbox{ in \cite{GMP}}.
\end{equation*}
However in the comparison results (Theorem \ref{uniqueness-thm-bd-domain} and Theorem \ref{uniqueness-thm2-bd-domain}) that we will prove in this paper the subsolutions and the supersolutions that we consider are $C^{2,1}(\Rr^n\times (0,T))$ functions and the condition $u_t\in L^1_{loc}(\Rr^n\times (0,T))$ is automatically satisfied. 

Asymptotic behaviour of the solution of 
\begin{equation*}\label{dp}
\left\{\begin{aligned}
u_t=&\Delta u^m\quad\,\mbox{ in }\Omega\times (0,T)\\
u=&0\qquad\,\,\,\mbox{ on }\partial\Omega\times (0,T)\\
u(x,0)=&u_0(x)\qquad\mbox{in }\Omega
\end{aligned}\right.
\end{equation*}
which vanishes at time $T>0$ where $0<m<1$, $n\ge 3$ and $u_0\ge 0$ is a function  on $\Omega$ was studied by G.~Akagi \cite{A1}, \cite{A2}, J.G.~Berryman and C.J.~Holland \cite{BH}, B.~Choi , R.J.~Mccann and C.~Seis \cite{CMS}, etc. Let $\mu_0>0$, $\mu_0\le f_1\in C^1(\2{\Omega})$ and $\mu_0\le v_0\in C^2(\2{\Omega}\setminus\{a_1\})$ satisfies
\begin{equation*}
\lambda_1|x-a_1|^{-\gamma_1}\le v_0(x)\le\lambda_1'|x-a_1|^{-\gamma_1'}\quad\forall\Omega\setminus\{a_1\}
\end{equation*} 
for some constants $\lambda_1'\ge\lambda_1>0$, $\gamma_1'\ge\gamma_1>\frac{2}{1-m}$.  When $n\ge 3$ and $0<m\le\frac{n-2}{n}$, existence and asymptotic large time behaviour of the Dirichlet blow-up solution of
\begin{equation*}
\left\{\begin{aligned}
u_t=&\Delta u^m\quad\,\mbox{ in }(\Omega\setminus\{a_1\})\times (0,\infty)\\
u=&f_1\qquad\,\mbox{ on }\partial\Omega\times (0,\infty)\\
u(a_1,t)=&\infty\qquad\,\forall t>0\\
u(x,0)=&v_0(x)\qquad\mbox{in }\Omega\setminus\{a_1\}
\end{aligned}\right.
\end{equation*}
has been proved by J.L.~Vazquez and M.~Winkler in \cite{VW1}, \cite{VW2}. 
When $n\ge 3$ and $0<m<\frac{n-2}{n}$, existence of finite blow-up points solution of
\begin{equation}\label{Dirichlet-blow-up-problem}
\left\{\begin{aligned}
u_t=&\Delta u^m\quad\,\mbox{ in }\widehat{\Omega}\times (0,\infty)\\
u=&f\qquad\,\,\,\mbox{ on }\partial\Omega\times (0,\infty)\\
u(a_i,t)=&\infty\qquad\,\forall t>0, i=1,2\dots,i_0\\
u(x,0)=&u_0(x)\qquad\mbox{in }\widehat{\Omega}
\end{aligned}\right.
\end{equation}
for any $0\le f\in L^{\infty}(\partial\Omega\times (0,\infty))$ and $0\le u_0\in L^p_{loc}(\2{\Omega}\setminus\{a_1,\dots, a_{i_0}\})$ for some constant $p>\frac{n(1-m)}{2}$ which satisfies
\begin{equation*}
u_0(x)\ge\lambda_i|x-a_i|^{-\gamma_i}\quad\forall |x-a_i|<\delta_1,i=1,\dots,i_0
\end{equation*} 
for some constants $0<\delta_1<\delta_0$, $\lambda_i>0$ and $\gamma_i>\frac{2}{1-m}$ for any $i=1,\dots, i_0$, has been proved by K.M.~Hui and S.~Kim \cite{HK2}. When $n\ge 3$, $0<m<\frac{n-2}{n}$ and $f$, $u_0$, also satisfy $f\ge\mu_0$ and $u_0\ge \mu_0$ for some constant $\mu_0>0$
and
\begin{equation*}
u_0(x)\le\lambda_i'|x-a_i|^{-\gamma_i'}\quad\forall |x-a_i|<\delta_1,i=1,\dots,i_0
\end{equation*} 
for some constants $\lambda_i'\ge\lambda_i>0$, $\gamma_i'\ge\gamma_i>\frac{2}{1-m}, i=1,\dots, i_0$, asymptotic large time behaviour of the finite blow-up points solution of \eqref{Dirichlet-blow-up-problem} has been proved by K.M.~Hui and S.~Kim in \cite{HK2} and \cite{H2}.

When $n\ge 3$ and $0<m\le\frac{n-2}{n}$, existence of finite blow-up points solutions of \eqref{fde} in bounded cylindrical domains was also proved by K.M.~Hui and Sunghoon Kim in \cite{HK1} using a different method when the initial value $u_0$ satisfies 
\begin{equation*}
u_0(x)\approx |x-a_i|^{-\gamma_i}\quad\mbox{ for }x\approx a_i\quad\forall i=1,2,\dots,i_0
\end{equation*} 
for some constants $\gamma_i>\max\left(\frac{n}{2m},\frac{n-2}{m}\right), i=1,2,\dots,i_0$.

\vspace{0.1in}

\noindent{\bf Outline of our results:} 

\begin{itemize}

\item We improve the existence theorems of \cite{HK2} (Theorem 1.1 and Theorem 1.2 of \cite{HK2}) to the existence of unique maximal solutions of \eqref{Dirichlet-blow-up-problem} (Theorem \ref{unique-soln-thm-bd-domain} and Theorem \ref{unique-soln-thm-bd-domain2}).

\item We extend the comparison theorems of \cite{H2} (Theorem 1.1 and Theorem 1.2 of \cite{H2}) by removing the requirement that the boundary values and the initial values must  be larger than some positive constant (Theorem \ref{uniqueness-thm-bd-domain} and Theorem \ref{uniqueness-thm2-bd-domain}). 

\item We extend the  asymptotic large time behaviour of the finite blow-up points solutions results of \cite{HK2} and \cite{H2} by removing the requirement that the boundary value $f$ and the initial value $u_0$ must  be larger than some positive constant (Theorem \ref{convergence-thm-bded-domain} and Theorem \ref{convergence-thm2-bded-domain}).
More precisely we prove the asymptotic large time behaviour of the finite blow-up points solution of \eqref{Dirichlet-blow-up-problem} (Theorem \ref{convergence-thm-bded-domain} and Theorem \ref{convergence-thm2-bded-domain}) for any $0\le f\in L^{\infty}(\partial\Omega\times (0,\infty))$ and $0\le u_0\in L^p_{loc}(\2{\Omega}\setminus\{a_1,\dots, a_{i_0}\})$ for some constant $p>\frac{n(1-m)}{2}$ which satisfies
\begin{equation}\label{initial-value-lower-upper-bd}
\frac{\lambda_i}{|x-a_i|^{\gamma_i}}\le u_0(x)\le \frac{\lambda_i'}{|x-a_i|^{\gamma_i'}} \qquad \forall 0<|x-a_i|<\delta_1,i=1, \cdots, i_0
\end{equation} 
for some constants $0<\delta_1<\delta_0$,  $\lambda_1$, $\cdots$, $\lambda_{i_0}$, $\lambda_1'$, $\cdots$, $\lambda_{i_0}'\in\mathbb{R}^+$ and 
\begin{equation}\label{gamma-gamma'-lower-upper-bd}
\frac{2}{1-m}<\gamma_i\le\gamma_i'<\frac{n-2}{m}\quad\forall i=1,2,\dots,i_0.
\end{equation}  

\item In the paper \cite{H2} K.M.~Hui constructed a solution of \eqref{Dirichlet-blow-up-problem} which oscillates between some fixed positive number and infinity as $t\to\infty$. A natural question to ask is whether there exist solutions of \eqref{Dirichlet-blow-up-problem} which oscillate between some functions on $\Omega$. We answer this question in the affirmative.  We will construct (Theorem \ref{oscillation-thm-bded-domain}) a solution of \eqref{Dirichlet-blow-up-problem} with appropriate lateral boundary value such that the solution of \eqref{Dirichlet-blow-up-problem} will oscillate between two given harmonic functions as $t\to\infty$.

\item We will prove the existence of minimal finite blow-up points solutions of \eqref{fde} in bounded cylindrical domains (Theorem \ref{unique-soln-thm3-bd-domain}) which also blow-up everywhere on the lateral boundary of the domain. Asymptotic large time behaviour of such solution is also prove in 
Theorem \ref{unique-soln-thm3-bd-domain}.
\end{itemize}

 More precisely we obtain the following results. The first four theorems are extensions of Theorem 2.3, Theorem 2.4 of \cite{H2} and Theorem 1.5 of \cite{HK2}.

\begin{thm}\label{unique-soln-thm-bd-domain}
Let $n\geq 3$, $0<m<\frac{n-2}{n}$, $0<\delta_1<\min (1,\delta_0)$,  $0\le f\in C^3(\partial\Omega\times (0,\infty))\cap L^{\infty}(\partial\Omega\times (0,\infty))$  and $0\le  u_0\in L_{loc}^p(\2{\Omega}\setminus\{a_1,\cdots,a_{i_0}\})$ for some constant  $p>\frac{n(1-m)}{2}$ be such that 
\eqref{initial-value-lower-upper-bd} holds for some constants  
\begin{equation}\label{gamma-gamma'-lower-bd}
\gamma_i'\ge\gamma_i>\frac{2}{1-m}\quad\forall i=1,2,\dots,i_0
\end{equation}  
 and $\lambda_1$, $\cdots$, $\lambda_{i_0}$, $\lambda_1'$, $\cdots$, $\lambda_{i_0}'\in\mathbb{R}^+$.  Let $\widehat{\Omega}$ be given by \eqref{omega-hat}. Then there exists a unique maximal solution $u$ of \eqref{Dirichlet-blow-up-problem} 
 such that for any constants $T>0$ and $\delta_2\in (0,\delta_1)$ there exist constants $C_1=C_1(T)>0$, $C_2=C_2(T)>0$, depending only on $\lambda_1$, 
$\cdots$, $\lambda_{i_0}$, $\lambda_1'$, $\cdots$, $\lambda_{i_0}'$, $\gamma_1$, $\cdots$, $\gamma_{i_0}$, 
$\gamma_1'$, $\cdots$, $\gamma_{i_0}'$, such that
\begin{equation}\label{u-lower-upper-blow-up-rate5}
\frac{C_1}{|x-a_i|^{\gamma_i}}\le u(x,t)\le\frac{C_2}{|x-a_i|^{\gamma_i'}} \quad \forall 0<|x-a_i|<\delta_2, 0<t<T, i=1,2,\dots,i_0
\end{equation}
holds. Moreover the following holds.

\begin{enumerate}[(i)]

\item If there exist constants $T_0'>T_0>0$ and $\mu_1>0$ such that
\begin{equation}\label{f-positive-after-t1}
f\ge\mu_1\quad\mbox{ on }\partial\Omega\times [T_0,T_0'),
\end{equation}
then for any $T_1\in (T_0,T_0')$ there exists a constant $\mu_2\in (0,\mu_1)$ such that
\begin{equation}\label{u-eventual-positive}
u(x,t)\ge\mu_2\quad\forall x\in\widehat{\Omega}, T_1\le t<T_0'.
\end{equation}

\item If there exists a constant $T_2\ge 0$ such that 
\begin{equation}\label{eq-condition-monotone-decreasingness-of-f2-34}
f(x,t)\mbox{ is monotone decreasing in $t$ on }\partial\Omega\times(T_2,\infty),
\end{equation} 
then $u$ satisfies
\begin{equation}\label{Aronson-Bernilan-ineqn}
u_t\leq\frac{u}{(1-m)(t-T_2)} \quad \mbox{ in }\widehat{\Omega}\times (T_2,\infty).
\end{equation}

\end{enumerate}

\end{thm}

\begin{thm}\label{unique-soln-thm-bd-domain2}
Let $n\ge 3$, $0<m<\frac{n-2}{n}$, $0<\delta_1<\min (1,\delta_0)$,  $0\le f\in  L^{\infty}(\partial\Omega\times (0,\infty))$  and $0\le  u_0\in L_{loc}^p(\2{\Omega}\setminus\{a_1,\cdots,a_{i_0}\})$ for some constant  $p>\frac{n(1-m)}{2}$ be such that 
\eqref{initial-value-lower-upper-bd} holds for some constants  $\gamma_i$, $\gamma_i'$, $i=1,\dots,i_0$, satisfying \eqref{gamma-gamma'-lower-upper-bd} and $\lambda_1$, $\cdots$, $\lambda_{i_0}$, $\lambda_1'$, $\cdots$, $\lambda_{i_0}'\in\mathbb{R}^+$.  Let $\widehat{\Omega}$ be given by \eqref{omega-hat}. Then there exists a unique maximal solution $u$ of \eqref{Dirichlet-blow-up-problem} 
 such that for any constants $T>0$ and $\delta_2\in (0,\delta_1)$ there exist constants $C_1=C_1(T)>0$, $C_2=C_2(T)>0$, depending only on $\lambda_1$, 
$\cdots$, $\lambda_{i_0}$, $\lambda_1'$, $\cdots$, $\lambda_{i_0}'$, $\gamma_1$, $\cdots$, $\gamma_{i_0}$, $\gamma_1'$, $\cdots$, $\gamma_{i_0}'$, such that
\eqref{u-lower-upper-blow-up-rate5} holds. Moreover the following holds.

\begin{enumerate}[(i)]

\item If there exists constants $T_0'>T_0>0$ and $\mu_1>0$ such that
\eqref{f-positive-after-t1} holds, then for any $T_1\in (T_0,T_0')$ there exists a constant $\mu_2\in (0,\mu_1)$ such that \eqref{u-eventual-positive} holds.

\item If there exists a constant $T_2\ge 0$ such that 
\eqref{eq-condition-monotone-decreasingness-of-f2-34} holds,
then $u$ satisfies \eqref{Aronson-Bernilan-ineqn}.
\end{enumerate}

\end{thm}

\begin{thm}\label{convergence-thm-bded-domain}
Let $n\ge 3$ and $0<m<\frac{n-2}{n}$. Let $0\le g\in C^3(\partial\Omega)$  and $\phi$ be the solution of 
\begin{equation}\label{harmonic-eqn}
\begin{cases}
\begin{aligned}
\Delta\phi&=0 \quad\,\,\,\, \mbox{ in }\Omega\\
\phi&=g^m \quad \mbox{ on }\partial\Omega.
\end{aligned}
\end{cases}
\end{equation}
Let $0\le u_0\in L^{p}_{loc}\left(\overline{\Omega}\setminus\left\{a_1,\cdots,a_{i_0}\right\}\right)$ for some constant $p>\frac{n(1-m)}{2}$ be
such that \eqref{initial-value-lower-upper-bd} holds for some constants $\gamma_i$, $\gamma_i'$, $i=1,\dots,i_0$, satisfying 
\eqref{gamma-gamma'-lower-upper-bd} and $0<\delta_1<\min (1,\delta_0)$, $\lambda_1$, $\cdots$, $\lambda_{i_0}$, $\lambda_1'$, $\cdots$, $\lambda_{i_0}'\in\mathbb{R}^+$. Let $f\in C^3(\partial\Omega\times\left(0,\infty\right)))\cap L^{\infty}(\partial\Omega\times\left(0,\infty\right)))$  be such that
\begin{equation}\label{f-t-infty-limit}
f\to g\quad\mbox{ uniformly in }C^3\left(\partial\Omega\right)\quad\mbox{ as }t\to\infty.
\end{equation}
 Let $\widehat{\Omega}$ be given by \eqref{omega-hat}.
Let $u$ be the unique maximal solution of \eqref{Dirichlet-blow-up-problem} given by Theorem \ref{unique-soln-thm-bd-domain}.  Then the following holds.
\begin{enumerate}[(i)]

\item If $g>0$ on $\partial\Omega$, then 
\begin{equation}\label{u-phi-1/m-limit1}
u(x,t)\to \phi^{\frac{1}{m}} \quad \mbox{ uniformly in } C^2(K) \quad \mbox{ as }t\to\infty
\end{equation}
holds  for any compact subset $K$ of $\overline{\Omega}\setminus\{a_1,\dots,a_{i_0}\}$.

\item If $g\not\equiv 0$ on $\partial\Omega$, 
\begin{equation}\label{f>g}
f\ge g\quad\mbox{ on }\partial\Omega\times (0,\infty)
\end{equation}
and
\begin{equation}\label{u0-lower-bd3}
u_0\ge\phi^{\frac{1}{m}}\quad\mbox{ on }\widehat{\Omega}
\end{equation}
holds, then \eqref{u-phi-1/m-limit1} holds for any compact subset $K$ of $\widehat{\Omega}$.

\item If $g\equiv 0$ on $\partial\Omega$, then

\begin{equation}\label{u-phi-1/m-limit2}
u(x,t)\to 0 \quad \mbox{ uniformly in }K \quad \mbox{ as }t\to\infty
\end{equation}
for any compact subset $K$ of $\widehat{\Omega}$.
\end{enumerate}
\end{thm}

\begin{thm}\label{convergence-thm2-bded-domain}
Let $n\ge 3$ and $0<m<\frac{n-2}{n}$. Let $0\le g\in C(\partial\Omega)$  and $\phi$ be the solution of 
\eqref{harmonic-eqn}.
Let $0\le u_0\in L^{p}_{loc}\left(\overline{\Omega}\setminus\left\{a_1,\cdots,a_{i_0}\right\}\right)$ for some constant $p>\frac{n(1-m)}{2}$ be
such that \eqref{initial-value-lower-upper-bd} holds for some constants $\gamma_i$, $\gamma_i'$, $i=1,\dots,i_0$, satisfying 
\eqref{gamma-gamma'-lower-upper-bd} and $0<\delta_1<\min (1,\delta_0)$, $\lambda_1$, $\cdots$, $\lambda_{i_0}$, $\lambda_1'$, $\cdots$, $\lambda_{i_0}'\in\mathbb{R}^+$. Let $f\in L^{\infty}(\partial\Omega\times\left(0,\infty\right)))$  be such that
\begin{equation*}
f\to g\quad\mbox{ uniformly in }L^{\infty}\left(\partial\Omega\right)\quad\mbox{ as }t\to\infty.
\end{equation*}
Let $\widehat{\Omega}$ be given by \eqref{omega-hat}. Let $u$ be the unique maximal solution of \eqref{Dirichlet-blow-up-problem} given by Theorem \ref{unique-soln-thm-bd-domain2}.  Then the following holds.

\begin{enumerate}[(i)]

\item If $g>0$ on $\partial\Omega$, then \eqref{u-phi-1/m-limit1}
holds  for any compact subset $K$ of $\widehat{\Omega}$.

\item If $g\not\equiv 0$ on $\partial\Omega$ and both \eqref{f>g} and \eqref{u0-lower-bd3}
holds, then \eqref{u-phi-1/m-limit1} holds for any compact subset $K$ of $\widehat{\Omega}$.

\item If $g\equiv 0$ on $\partial\Omega$, then \eqref{u-phi-1/m-limit2} holds
for any compact subset $K$ of $\widehat{\Omega}$.

\end{enumerate}

\end{thm}

\begin{thm}\label{oscillation-thm-bded-domain}
Let $n\ge 3$ and $0<m<\frac{n-2}{n}$. Let $g_1,g_2\in C(\partial\Omega)$, $g_2>0$,  $g_1>0$,  and $\phi_1$, $\phi_2$, be the solutions of \eqref{harmonic-eqn} with $g=g_1, g_2$, respectively. Let
\begin{equation*}
0<\mu_0<\min \left(\min_{\partial\Omega}g_1,\min_{\partial\Omega}g_2\right)
\end{equation*}
be a constant.
Let $0\le u_0\in L^{p}_{loc}\left(\overline{\Omega}\setminus\left\{a_1,\cdots,a_{i_0}\right\}\right)$ for some constant $p>\frac{n(1-m)}{2}$ be such that \eqref{initial-value-lower-upper-bd} holds for some constants $\gamma_i$, $\gamma_i'$, $i=1,\dots,i_0$, satisfying 
\eqref{gamma-gamma'-lower-upper-bd} and $0<\delta_1<\min (1,\delta_0)$, $\lambda_1$, $\cdots$, $\lambda_{i_0}$, $\lambda_1'$, $\cdots$, $\lambda_{i_0}'\in\mathbb{R}^+$. Let $\widehat{\Omega}$ be given by \eqref{omega-hat}. Then there exist a function $f\in L^{\infty}(\partial\Omega\times\left(0,\infty\right)))$ and an increasing sequence $\{t_i\}_{i=1}^{\infty}$, $t_i\to\infty$ as $i\to\infty$, such that if $u$ is the maximal solution  of \eqref{Dirichlet-blow-up-problem} given by Theorem \ref{unique-soln-thm-bd-domain2}, then

\begin{equation*}
\left\{\begin{aligned}
&u(x,t_{2i-1})\to \phi_1^{\frac{1}{m}} \quad \mbox{ in } C^2(K) \quad \mbox{ as }i\to\infty\\
&u(x,t_{2i})\to \phi_2^{\frac{1}{m}} \qquad \mbox{in } C^2(K) \quad \mbox{ as }i\to\infty
\end{aligned}\right.
\end{equation*}
for any compact subset $K$ of $\widehat{\Omega}$.

\end{thm}

\begin{thm}\label{unique-soln-thm3-bd-domain}
Let $n\geq 3$, $0<m<\frac{n-2}{n}$, $0<\delta_1<\min (1,\delta_0)$  and $0\le  u_0\in L_{loc}^p(\2{\Omega}\setminus\{a_1,\cdots,a_{i_0}\})$ for some constant  $p>\frac{n(1-m)}{2}$ be such that 
\eqref{initial-value-lower-upper-bd} holds for some constants $\gamma_i$, $\gamma_i'$, $i=1,\dots,i_0$, satisfying \eqref{gamma-gamma'-lower-bd}
 and $\lambda_1$, $\cdots$, $\lambda_{i_0}$, $\lambda_1'$, $\cdots$, $\lambda_{i_0}'\in\mathbb{R}^+$. Then there exists a unique minimal  solution $u$ of 
\begin{equation}\label{Dirichlet-blow-up-origin-bdary-problem}
\left\{\begin{aligned}
u_t=&\Delta u^m\quad\,\mbox{ in }\widehat{\Omega}\times (0,\infty)\\
u=&\infty\qquad\mbox{ on }\partial\Omega\times (0,\infty)\\
u(a_i,t)=&\infty\qquad\forall t>0, i=1,2\dots,i_0\\
u(x,0)=&u_0(x)\qquad\mbox{in }\widehat{\Omega}
\end{aligned}\right.
\end{equation}
such that for any constants $T>0$ and $\delta_2\in (0,\delta_1)$ there exist constants $C_1=C_1(T)>0$, $C_2=C_2(T)>0$, depending only on $\lambda_1$, 
$\cdots$, $\lambda_{i_0}$, $\lambda_1'$, $\cdots$, $\lambda_{i_0}'$, $\gamma_1$, $\cdots$, $\gamma_{i_0}$, $\gamma_1'$, $\cdots$, $\gamma_{i_0}'$, such that \eqref{u-lower-upper-blow-up-rate5} holds.
Moreover  $u$ satisfies \eqref{Aronson-Bernilan-ineqn} with $T_2=0$ and
\begin{equation}\label{u-t-inftyu-infty}
u(x,t)\to\infty\quad\mbox{ uniformly on }\Omega_{\delta}\quad\mbox{ as }t\to\infty\quad\forall 0<\delta<\delta_0.
\end{equation}

\end{thm}

\begin{rmk}
The integrability condition $0\le  u_0\in L_{loc}^p(\2{\Omega}\setminus\{a_1,\cdots,a_{i_0}\})$ for some constant  $p>\frac{n(1-m)}{2}$ is necessary since this condition together with $f\in L^{\infty}(\1\Omega\times (0,\infty))$ implies that the solution $u$ of  \eqref{Dirichlet-blow-up-problem}  locally satisfies a $L^{\infty}-L^p$ regularizing result in terms of the local $L^p$ norm of the initial value $u_0$ and $L^{\infty}$ norm of $f$ (Lemma \ref{lem-Cor-2-2-of-Hs1} and Lemma \ref{L-infinity-initial-Lp-bd1}).
\end{rmk}

\begin{rmk}
In the proof of Theorem \ref{unique-soln-thm-bd-domain} and Theorem \ref{unique-soln-thm-bd-domain2} we will construct the  solution of  \eqref{Dirichlet-blow-up-problem} as the limit of a monotone decreasing sequence of solutions  of
\begin{equation}\label{blow-up-problem}
\left\{\begin{aligned}
u_t=&\Delta u^m\quad\,\mbox{ in }\widehat{\Omega}\times (0,\infty)\\
u(a_i,t)=&\infty\qquad\,\forall t>0, i=1,2\dots,i_0\\
\end{aligned}\right.
\end{equation}
which have initial values strictly greater than $u_0$ and lateral boundary values strictly greater than $f$. Since there is comparison results (Theorem \ref{uniqueness-thm-bd-domain} and Theorem \ref{uniqueness-thm2-bd-domain}) between any solution of  \eqref{Dirichlet-blow-up-problem} and  this monotone decreasing sequence of solutions of \eqref{blow-up-problem}. Hence this constructed  solution of  \eqref{Dirichlet-blow-up-problem}  must be  maximal solution of  \eqref{Dirichlet-blow-up-problem} by comparison argument. On the other hand if we construct  solution of  \eqref{Dirichlet-blow-up-problem} as the limit of a monotone increasing sequence of solutions of \eqref{blow-up-problem} which have initial values less than $u_0$ and lateral boundary values less than $f$. Since there is no comparison result between any solution of  \eqref{Dirichlet-blow-up-problem} and  this monotone increasing sequence of solutions of \eqref{blow-up-problem}. Hence it is not clear whether minimal solution of  \eqref{Dirichlet-blow-up-problem} exists.
\end{rmk}

The plan of the paper is as follows. For the readers' convenience in section 2 we recall some results of \cite{H1}, \cite{H2} and \cite{HK2} that is cited in this paper.
In section 3 we will prove Theorem  \ref{unique-soln-thm-bd-domain} and Theorem  \ref{unique-soln-thm-bd-domain2}. We will prove Theorem \ref{convergence-thm-bded-domain} and Theorem \ref{convergence-thm2-bded-domain} in section 4. We will prove Theorem \ref{oscillation-thm-bded-domain} and Theorem \ref{unique-soln-thm3-bd-domain} in section 5. Unless stated otherwise we will assume that  $n\ge 3$ and $0<m<\frac{n-2}{n}$ for the rest of the paper.

We start with some definitions. For any $a\in\mathbb{R}^+$, let $a_+=\max (a,0)$. 

\begin{defn}
For any $t_2>t_1$, we say that $u$ is a solution  (subsolution, supersolution respectively) of \eqref{fde} in $\Omega\times(t_1,t_2)$ if $u\in C^{2,1}(\Omega\times (t_1,t_2))$ is positive in $\Omega\times(t_1,t_2)$  and  satisfies
\begin{equation*}
u_t=\Delta u^m \quad\mbox{ in }\Omega\times (t_1,t_2)\quad (\le, \ge, \,\mbox{ respectively}).
\end{equation*}
\end{defn}

\begin{defn}\label{dirichlet-problem-soln-defn}
For any $0\leq f\in L^{\infty}(\partial\Omega\times(0,T))$ and $0\le u_0\in L^{1}_{loc}(\Omega)$, we say that $u$ is a solution (subsolution, supersolution respectively) of 
\begin{equation}\label{fde-Dirichlet-problem}
\begin{cases}
\begin{aligned}
u_t&=\Delta u^m \qquad \mbox{in $\Omega\times(0,T)$}\\
u&=f \qquad \quad \mbox{on $\partial\Omega\times(0,T)$}\\
u(x,0)&=u_0(x) \quad\,\,\, \mbox{in $\Omega$}.
\end{aligned}
\end{cases}
\end{equation}
if $u$ is a solution  (subsolution, supersolution respectively) of \eqref{fde} in $\Omega\times(0,T)$ which satisfies 
\begin{equation*}
\|u(\cdot,t)-u_0\|_{L^1(\Omega)}\to 0\quad \mbox{ as }t\to 0 
\end{equation*} 
and the boundary condition is satisfied in the sense that
\begin{equation*}\label{eq-alinged-formula-for-weak-solution-of-u}
\int_{t_1}^{t_2}\int_{\Omega}(u\eta_t+u^m\Delta\eta)\,dx\,dt=\int_{t_1}^{t_2}\int_{\partial\Omega}f^m\frac{\partial\eta}{\partial\nu}\,d\sigma\, dt+\int_{\Omega}u\eta\,dx\bigg|_{t_1}^{t_2}
\end{equation*}
($\geq, \leq$ respectively) holds for any $0<t_1<t_2<T$ and $\eta\in C_c^2(\overline{\Omega}\times(0,T))$ satisfying $\eta=0$ on $\partial\Omega\times(0,T)$.
\end{defn}

\begin{defn}
For any $T>0$, $0\leq f\in L^{\infty}(\partial\Omega\times(0,T))$ and $0\le u_0\in L^1_{loc}(\widehat{\Omega})$ where $\widehat{\Omega}$ is given by \eqref{omega-hat}, we say that $u$ is a solution (subsolution, supersolution respectively) of 
\begin{equation}\label{puntured-Dirichlet-problem}
\begin{cases}
\begin{aligned}
u_t&=\Delta u^m \qquad \mbox{in $\widehat{\Omega}\times(0,T)$}\\
u(x,t)&=f \qquad \quad \mbox{on $\partial\Omega\times(0,T)$}\\
u(x,0)&=u_0(x) \quad\,\,\, \mbox{in $\widehat{\Omega}$}.
\end{aligned}
\end{cases}
\end{equation} 
if $u$ is a solution  (subsolution, supersolution respectively) of \eqref{fde} in $\widehat{\Omega}\times(0,T)$ which satisfies
\begin{equation}\label{u-initial-value}
\|u(\cdot,t)-u_0\|_{L^1(K)}\to 0\quad \mbox{ as }t\to 0 
\end{equation} 
for any compact set $K\subset\overline{\Omega}\setminus\left\{a_1,\cdots,a_{i_0}\right\}$  and
\begin{equation*}\label{eq-formula-for-weak-sols-of-main-problem}
\begin{aligned}
&\int_{t_1}^{t_2}\int_{\widehat{\Omega}}\left(u\eta_t+u^m\Delta\eta\right)\,dxdt\\
&\qquad \qquad =\int_{t_1}^{t_2}\int_{\partial\Omega}f^m\frac{\partial\eta}{\partial\nu}\,d\sigma dt+\int_{\widehat{\Omega}}u(x,t_2)\eta(x,t_2)\,dx-\int_{\widehat{\Omega}}u(x,t_1)\eta(x,t_1)\,dx
\end{aligned}
\end{equation*}
($\geq, \leq$ respectively) for any $0<t_1<t_2<T$ and $\eta\in C_c^2((\overline{\Omega}\setminus\left\{a_1,\cdots,a_{i_0}\right\})\times(0,T))$ satisfying $\eta\equiv 0$ on $\partial\Omega\times(0,T)$.  
\end{defn}

\begin{defn}
We say that $u$ is a solution (subsolution, supersolution respectively) of \eqref{Dirichlet-blow-up-problem} if $u$ is a solution (subsolution, supersolution respectively) of 
\eqref{puntured-Dirichlet-problem} and 
\begin{equation}\label{u-blow-up-at-ai}
u(x,t)\to\infty\quad\mbox{ as }x\to a_i\quad\forall t>0, i=1,\dots,i_0.
\end{equation}
\end{defn}

\begin{defn}
We say that $u$ is a maximal solution of \eqref{Dirichlet-blow-up-problem} 
if $u$ is a solution of \eqref{Dirichlet-blow-up-problem} and for any solution $v$ of \eqref{Dirichlet-blow-up-problem}, $v\le u$ in $\widehat{\Omega}\times (0,T)$.
\end{defn}

\begin{defn}
We say that $u$ is a solution of \eqref{Dirichlet-blow-up-origin-bdary-problem} if $u$ is a solution  of
\eqref{fde} in $\widehat{\Omega}\times (0,\infty)$  which satisfies  \eqref{u-initial-value} for any compact set $K\subset\widehat{\Omega}$, \eqref{u-blow-up-at-ai} and
\begin{equation*}
\underset{\substack{(y,s)\to (x,t)\\
(y,s)\in\Omega\times (0,\infty)}}{\lim}u(y,s)=\infty
\quad\forall (x,t)\in\partial\Omega\times (0,\infty).
\end{equation*}
\end{defn}

\begin{defn}
We say that $u$ is a minimal solution of \eqref{Dirichlet-blow-up-origin-bdary-problem} 
if $u$ is a solution of \eqref{Dirichlet-blow-up-origin-bdary-problem} and for any solution $v$ of \eqref{Dirichlet-blow-up-origin-bdary-problem}, $v\ge u$ in $\widehat{\Omega}\times (0,T)$.
\end{defn}

\section{Preliminaries}
\setcounter{equation}{0}
\setcounter{thm}{0}

In this section we recall some results of \cite{H1}, \cite{H2} and \cite{HK2} that are cited in this paper.

\begin{thm}[Theorem 1.1 of \cite{H2}]
Let $n\geq 3$, $0<m<\frac{n-2}{n}$, $0<\delta_1<\min (1,\delta_0)$, $\mu_0>0$, $f_1, f_2\in C^3(\partial\Omega\times (0,\infty))\cap L^{\infty}(\partial\Omega\times (0,\infty))$ be such that $f_2\ge f_1\ge\mu_0$ on $\partial\Omega\times (0,\infty)$ and 
\begin{equation}\label{u0-1-2-ineqn0}
\mu_0\leq u_{0,1}\le u_{0,2}\in L_{loc}^p(\2{\Omega}\setminus\{a_1,\cdots,a_{i_0}\})\quad\mbox{ for some constant  }p>\frac{n(1-m)}{2}
\end{equation}
be such that 
\begin{equation}\label{2initial-value-lower-upper-bd0}
\frac{\lambda_i}{|x-a_i|^{\gamma_i}}\le u_{0,1}(x)\le u_{0,2}\le \frac{\lambda_i'}{|x-a_i|^{\gamma_i'}} \qquad \forall 0<|x-a_i|<\delta_1,i=1, \cdots, i_0
\end{equation}
holds for some constants $\lambda_1$, $\cdots$, $\lambda_{i_0}$, $\lambda_1'$, $\cdots$, $\lambda_{i_0}'\in\Rr^+$ and 
\begin{equation*}
\gamma_i'\ge\gamma_i>\frac{2}{1-m}\quad\forall i=1,2,\dots,i_0.
\end{equation*}  
Suppose $u_1$, $u_2$, are the solutions of \eqref{Dirichlet-blow-up-problem} with $u_0=u_{0,1}, u_{0,2}$, $f=f_1, f_2$, respectively which satisfy
\begin{equation}\label{u1-2-lower-bd0}
u_j(x,t)\ge\mu_0\quad\forall x\in \widehat{\Omega}, t>0, j=1,2
\end{equation}
such that for any constants $T>0$ and $\delta_2\in (0,\delta_1)$ there exist constants $C_1=C_1(T)>0$, $C_2=C_2(T)>0$, such that
\begin{equation}\label{u-lower-upper-blow-up-rate0}
\frac{C_1}{|x-a_i|^{\gamma_i}}\le u_j(x,t)\le\frac{C_2}{|x-a_i|^{\gamma_i'}} 
\end{equation}
holds for any  $0<|x-a_i|<\delta_2$, $0<t<T$, $i=1,2,\dots,i_0, j=1,2$. Suppose $u_1$, $u_2$, also satisfy
\begin{equation}\label{initial-value-l1-sense0}
\|u_i(\cdot,t)-u_{0,i}\|_{L^1(\Omega_{\delta})}\to 0\quad\mbox{ as }t\to 0\quad\forall 0<\delta<\delta_0, i=1,2.
\end{equation}
Then 
\begin{equation*}
u_1(x,t)\le u_2(x,t)\quad\forall x\in \widehat{\Omega}, t>0.
\end{equation*}
\end{thm}

\begin{thm}[Theorem 1.2 of \cite{H2}]
Let $n\geq 3$, $0<m<\frac{n-2}{n}$, $0<\delta_1<\min (1,\delta_0)$, $\mu_0>0$, $\mu_0\le f_1\le f_2\in L^{\infty}(\partial\Omega\times (0,\infty))$ and \eqref{u0-1-2-ineqn0},  \eqref{2initial-value-lower-upper-bd0},
hold for some constants $\lambda_1$, $\cdots$, $\lambda_{i_0}$, $\lambda_1'$, $\cdots$, $\lambda_{i_0}'\in\Rr^+$ satisfying \eqref{gamma-gamma'-lower-upper-bd}.
Suppose $u_1$, $u_2$, are the solutions of \eqref{Dirichlet-blow-up-problem} with $u_0=u_{0,1}, u_{0,2}$, $f=f_1, f_2$, respectively which 
satisfy  \eqref{u1-2-lower-bd0} and \eqref{initial-value-l1-sense0} 
such that for any constants $T>0$ and $\delta_2\in (0,\delta_1)$ there exist constants $C_1=C_1(T)>0$, $C_2=C_2(T)>0$, such that \eqref{u-lower-upper-blow-up-rate0} holds.
Then 
\begin{equation*}
u_1(x,t)\le u_2(x,t)\quad\forall x\in \widehat{\Omega}, t>0.
\end{equation*}
\end{thm}

For any $m\in\Rr$, we let $\phi_m(u)=u^m/m$ if $m\ne 0$ and $\phi_m(u)=\log u$ if $m=0$. 

\begin{lem}[Lemma 1.7 of \cite{H1}]
Let $m_0<0<\3_1<1$ and $m\in [m_0,1-\3_1]$. Suppose $u$ is a solution of 
\begin{equation*}
u_t=\Delta \phi_m(u)
\end{equation*}
in $\Omega\times (0,T)$ with initial value $0\le u_0\in L_{loc}^p(\Omega)$ for some constant 
\begin{equation*}
p>\max\left(1,(1-m_0)\max (1,n/2)\right).
\end{equation*}
Then for any $B_{R_1}(x_0)\subset\2{B_{R_2}(x_0)}\subset\Omega$ there exists a constant $C>0$ 
such that
$$
\int_{B_{R_1}(x_0)}u(x,t)^pdx
\le C\biggl\{t^{p/(1-m_0)}+t^{p/\3_1}+\int_{B_{R_2}(x_0)}
u_0^pdx\biggr\}
$$
holds for any $0\le t<T$, $m\in [m_0,1-\3_1]$.
\end{lem}

\begin{thm}[Theorem 1.1 of \cite{HK2}]\label{first-main-existence-thm}
Let $n\geq 3$, $0<m<\frac{n-2}{n}$, $0<\delta_1<\delta_0$, $0\le f\in L^{\infty}(\partial\Omega\times[0,\infty))$ and $0\leq u_0\in L_{loc}^p(\2{\Omega}\setminus\{a_1,\cdots,a_{i_0}\})$ for some constant  $p>\frac{n(1-m)}{2}$ be such that \begin{equation*}\label{u0-lower-blow-up-rate}
u_0(x)\geq \frac{\lambda_i}{|x-a_i|^{\gamma_i}} \qquad \forall 0<|x-a_i|<\delta_1,\,\,i=1, \cdots, i_0
\end{equation*}
holds for some constants $\lambda_1$, $\cdots$, $\lambda_{i_0}\in\Rr^+$ and $\gamma_1,\cdots,\gamma_{i_0}\in (\frac{2}{1-m},\infty)$. Then there exists a solution $u$ of 
\eqref{Dirichlet-blow-up-problem}
such that for any $T>0$ and $\delta_2\in (0,\delta_1)$ there exists a constant $C_1>0$  such that
\begin{equation*}\label{eq-lower-limit-of-u-near-blow-up-points}
u(x,t) \geq \frac{C_1}{|x-a_i|^{\gamma_i}}\quad \forall 0<\left|x-a_i\right|<\delta_2, 0<t<T.
\end{equation*}
Moreover if there exists a constant $T_2\ge 0$ such that \eqref{eq-condition-monotone-decreasingness-of-f2-34} holds, 
then $u$ satisfies \eqref{Aronson-Bernilan-ineqn}.
\end{thm}

\begin{thm}[Theorem 1.5 of \cite{HK2}]
Suppose that $n\geq 3$, $0<m<\frac{n-2}{n}$ and $\mu_0>0$. Let $\mu_0\leq u_0\in L^{p}_{loc}\left(\overline{\Omega}\backslash\left\{a_1,\cdots,a_{i_0}\right\}\right)$ for some constant $p>\frac{n(1-m)}{2}$ such that  \eqref{initial-value-lower-upper-bd} holds for some constants $\lambda_1$, $\cdots$, $\lambda_{i_0}$, $\lambda_1'$, $\cdots$, $\lambda_{i_0}'\in\mathbb{R}^+$ and $\gamma_i$, $\gamma_i'$, $i=1,\dots,i_0$, satisfying \eqref{gamma-gamma'-lower-upper-bd}.  Let $f\in L^{\infty}(\partial\Omega\times\left(0,\infty\right))\cap C^3(\partial\Omega\times\left(T_1,\infty\right))$ for some constant $T_1>0$ satisfy 
\begin{equation*}
f\geq\mu_0 \quad \mbox{ on }\partial\Omega\times\left(0,\infty\right)
\end{equation*}
and \eqref{f-t-infty-limit} for some function $g\in C^{3}\left(\partial\Omega\right)$, $g\geq\mu_0$ on $\partial\Omega$. Let $u$ be the solution of \eqref{Dirichlet-blow-up-problem} given by Theorem \ref{first-main-existence-thm}. Let $\psi$ be the solution of \eqref{harmonic-eqn}. Then \eqref{u-phi-1/m-limit1} holds
for any compact subset $K$ of $\overline{\Omega}\backslash\left\{a_1,\cdots,a_{i_0}\right\}$.
\end{thm}

\begin{lem}[Lemma 2.9 of \cite{HK2}]
Let $n\ge 3$, $0<m\leq\frac{n-2}{n}$, $0\le f\in L^{\infty}(\partial\Omega\times[0,\infty))$  and $0\le u_0\in L^p_{loc}(\2{\Omega}\setminus\{a_1,\cdots,a_{i_0}\})$ for some constant $p>\frac{n(1-m)}{2}$. Suppose $u$ is a solution of  \eqref{puntured-Dirichlet-problem}.
Then for any $0<\delta_6<\delta_5<\delta_0$ and $0<t_1<T$ there exist constants $C>0$ and $\theta>0$ such that
\begin{equation*}
\left\|u\right\|_{L^{\infty}(\Omega_{\delta_5}\times[t_1,T))}\leq C\left(k_f^p|\Omega|+\int_{\Omega_{\delta_6}}u_0^p\,dx\right)^{\theta/p}+k_f
\end{equation*}
where $k_f=\max (1,\|f\|_{L^{\infty}})$.
\end{lem}

\begin{lem}[Lemma 3.2 of \cite{HK2}]\label{lem-local-upper-bound-of-u-near-singular-points}
Let $n\geq 3$, $0<m<\frac{n-2}{n}$, $0<\delta_1<\min(1,\delta_0)$, $0\le f\in L^{\infty}(\partial\Omega\times[0,\infty))$ and $0\leq u_0\in L^p_{loc}(\2{\Omega}\setminus\{a_1,\cdots,a_{i_0}\})$ for some constant $p>\frac{n(1-m)}{2}$ such that  \eqref{initial-value-lower-upper-bd} holds for some constants $\lambda_1$, $\cdots$, $\lambda_{i_0}$, $\lambda_1'$, $\cdots$, $\lambda_{i_0}'\in\mathbb{R}^+$ and $\gamma_i$, $\gamma_i'$, $i=1,\dots,i_0$, satisfying \eqref{gamma-gamma'-lower-upper-bd}.
Let $u$ be the solution of \eqref{Dirichlet-blow-up-problem} given by Theorem \ref{first-main-existence-thm}. Then for any $0<\delta_2<\delta_0$ and $t_0>0$ there exist constants $C_2>0$ and $C_3>0$ such that
\begin{equation*}
u(x,t)\le C_2 \quad \forall x\in\overline{\Omega_{\delta_2}}\times[t_0,\infty)
\end{equation*} 
and
\begin{equation*}
u(x,t)\le C_3|x-a_i|^{-\gamma_i'}\quad\forall 0<|x-a_i|\leq\delta_2,t\ge t_0,i=1,\cdots,i_0
\end{equation*}
hold.
\end{lem}

\begin{rmk}[Remark 3.7 of \cite{HK2}]
If $f\in L^{\infty}(\partial\Omega\times(0,\infty))$, $g\in C(\1\Omega)$ and 
\begin{equation*}
f(x.t)\to g(x) \quad\mbox{ uniformly in }L^{\infty}(\partial\Omega)\mbox{ as }t\to\infty,
\end{equation*}
then  the solution $u$ of \eqref{Dirichlet-blow-up-problem} given by Theorem 1.1 of \cite{HK2} satisfy \eqref{u-phi-1/m-limit1} for any compact set $K\subset\widehat{\Omega}$. Moreover \begin{equation*}
u(x,t)\to \psi^{\frac{1}{m}} \quad \mbox{ in }L_{loc}^{\infty}(\overline{\Omega}\backslash\left\{a_1,\cdots,a_{i_0}\right\}) \quad \mbox{ as }t\to\infty.
\end{equation*} 
\end{rmk}

\section{Existence of maximal blow-up solutions}
\setcounter{equation}{0}
\setcounter{thm}{0}

In this section we will use a modification of the argument of \cite{HK2} and \cite{H2} to prove  the existence of maximal solution of \eqref{Dirichlet-blow-up-problem}. We first extend Theorem 1.1 and Theorem 1.2 of \cite{H2}. We start with a technical lemma.

\begin{lem}\label{up-u0-p-integral-ineqn-lem}
Let $n\ge 3$, $0<m<1$, $p>\frac{n(1-m)}{2}$, $0\le f\in L^{\infty}(\partial\Omega\times (0,T))$  and $0\le u_0\in L^{\infty}(\Omega)$. Suppose $u\in L^{\infty}(\Omega\times (0,T))$ is a solution of  \eqref{fde-Dirichlet-problem}.
Then for any $0<\delta'<\delta<\delta_0$ there exists a  constant $C>0$ depending only on $p$, $m$, $\delta$ and $\delta'$ such that
\begin{equation}\label{up-u0-p-integral-ineqn}
\int_{\Omega_{\delta}}u(x,t)^p\,dx\le C\left(\int_{\Omega_{\delta'}}u_0^p\,dx+t^{\frac{p}{1-m}}+ \|f\|_{L^{\infty}(\1\Omega\times (0,T))}^p\right)\quad\forall 0<t<T
\end{equation}
where $\Omega_{\delta}$, $\Omega_{\delta'}$, is given by \eqref{d-delta}.
\end{lem}
\begin{proof}
We will use a modification of the proof of Lemma 1.7 of \cite{H1} to prove this lemma. Let $0\le\phi_1\in C^{\infty}_0(\overline{\Omega}\setminus\{a_1,\cdots,a_{i_0}\})$, $0\le\phi_1\le 1$, be such that $\phi_1(x)=1$ for any $x\in \overline{\Omega_{\delta}}$ and $\phi_1(x)=0$ for any $x\in \Omega\setminus\Omega_{\delta'}$. Let $\phi_2=\phi_1^{\alpha}$ for some constant $\alpha>\frac{2p}{1-m}$ and $k>\|f\|_{L^{\infty}}$. Let $\widehat{\Omega}$ be given by \eqref{omega-hat}. Then
\begin{align}\label{up-integral-ineqn0}
&\frac{\partial}{\partial t}\left(\int_{\widehat{\Omega}}(u-k)_+^p\phi_2\,dx\right)\notag\\
=&p\int_{\widehat{\Omega}}(u-k)_+^{p-1}u_t\phi_2\,dx\notag\\
=&p\int_{\widehat{\Omega}}(u-k)_+^{p-1}\phi_2\Delta u^m\,dx\notag\\
=&-p\int_{\widehat{\Omega}}\nabla u^m\cdot\nabla[(u-k)_+^{p-1}\phi_2]\,dx\notag\\
=&-pm\left\{(p-1)\int_{\widehat{\Omega}}u^{m-1}(u-k)_+^{p-2}|\nabla u|^2\phi_2\,dx
+\int_{\widehat{\Omega}}u^{m-1}(u-k)_+^{p-1}\nabla u\cdot\nabla\phi_2\,dx\right\}.
\end{align}
Since
\begin{align*}
\left|\int_{\widehat{\Omega}}u^{m-1}(u-k)_+^{p-1}\nabla u\cdot\nabla\phi_2\,dx\right|
\le&(p-1)\int_{\widehat{\Omega}}u^{m-1}(u-k)_+^{p-2}|\nabla u|^2\phi_2\,dx\notag\\
&\qquad +\frac{1}{4(p-1)}\int_{\widehat{\Omega}}u^{m-1}(u-k)_+^p|\nabla\phi_2|^2\phi_2^{-1}\,dx,
\end{align*}
by \eqref{up-integral-ineqn0} and  H\"older's inequality with exponents $\frac{p}{1-m}$ and $\frac{p}{p+m-1}$,
\begin{align}\label{up-integral-ineqn}
&\frac{\partial}{\partial t}\left(\int_{\widehat{\Omega}}(u-k)_+^p\phi_2\,dx\right)\notag\\
\le&\frac{pm}{4(p-1)}\int_{\widehat{\Omega}}u^{m-1}(u-k)_+^p|\nabla\phi_2|^2\phi_2^{-1}\,dx\notag\\
\le&\frac{pm}{4(p-1)}\int_{\widehat{\Omega}}[(u-k)_+^p\phi_2]^{\frac{p+m-1}{p}}|\nabla\phi_2|^2
\phi_2^{\frac{1-m}{p}-2}\,dx\notag\\
\le&\frac{pm}{4(p-1)}\left(\int_{\widehat{\Omega}}\left(|\nabla\phi_2|^2\phi_2^{\frac{1-m}{p}-2}\right)^{\frac{p}{1-m}}\,dx\right)^{\frac{1-m}{p}}\left(\int_{\widehat{\Omega}}(u-k)_+^p\phi_2\,dx\right)^{1-\frac{1-m}{p}}\quad\forall 0<t<T.
\end{align}
Since
\begin{equation*}
|\nabla\phi_2|^2\phi_2^{\frac{1-m}{p}-2}=\alpha^2\phi_1^{\frac{(1-m)\alpha}{p}-2}|\nabla\phi_1|^2
\le \alpha^2|\nabla\phi_1|^2,
\end{equation*}
we have
\begin{equation*}
\int_{\widehat{\Omega}}\left(|\nabla\phi_2|^2\phi_2^{\frac{1-m}{p}-2}\right)^{\frac{p}{1-m}}\,dx
\le\alpha^{\frac{2p}{1-m}}\int_{\widehat{\Omega}}|\nabla\phi_1|^{\frac{2p}{1-m}}\,dx<\infty.
\end{equation*}
Then, by \eqref{up-integral-ineqn},
\begin{align}\label{u-k-integral-ineqn2}
&\frac{\partial}{\partial t}\left(\int_{\widehat{\Omega}}(u-k)_+^p\phi_2\,dx\right)
\le C\left(\int_{\widehat{\Omega}}(u-k)_+^p\phi_2\,dx\right)^{1-\frac{1-m}{p}}\notag\\
\Rightarrow\quad&\left(\int_{\widehat{\Omega}}(u-k)_+^p\phi_2\,dx\right)^{\frac{1-m}{p}-1}\frac{\partial}{\partial t}\left(\int_{\widehat{\Omega}}(u-k)_+^p\phi_2\,dx\right)\le C.
\end{align}
where $\widehat{\Omega}$ is given by \eqref{omega-hat}. 
Integrating \eqref{u-k-integral-ineqn2} over $(t_1,t)$, $0<t_1<t<T$,
\begin{equation*}
\int_{\widehat{\Omega}}(u(x,t)-k)_+^p\phi_2\,dx\le\left\{Ct+ \left(\int_{\widehat{\Omega}}(u(x,t_1)-k)_+^p\phi_2\,dx\right)^{\frac{1-m}{p}}\right\}^{\frac{p}{1-m}}.
\end{equation*}
Hence
\begin{align*}
\int_{\Omega_{\delta}}(u(x,t)-k)_+^p\,dx\le&\left\{Ct+ \left(\int_{\Omega_{\delta'}}(u(x,t_1)-k)_+^p\,dx\right)^{\frac{1-m}{p}}\right\}^{\frac{p}{1-m}}\notag\\
\le&C'\left(\int_{\Omega_{\delta'}}u(x,t_1)^p\,dx+t^{\frac{p}{1-m}}+ k^p\right)
\end{align*}
holds for any $0<t_1\le t<T$, $k>\|f\|_{L^{\infty}(\1\Omega\times (0,T))}$. Thus
\begin{align}\label{u-k-integral-ineqn}
\int_{\Omega_{\delta}}u(x,t)^p\,dx
\le&C\left(\int_{\Omega_{\delta'}}u(x,t_1)^p\,dx+t^{\frac{p}{1-m}}+ k^p\right)
\end{align}
holds for any $0<t_1\le t<T$, $k>\|f\|_{L^{\infty}(\1\Omega\times (0,T))}$.
Letting $k\searrow\|f\|_{L^{\infty}(\1\Omega\times (0,T))}$ in \eqref{u-k-integral-ineqn},
\begin{equation}\label{up-fp-integral-ineqn}
\int_{\Omega_{\delta}}u(x,t)^p\,dx
\le C\left(\int_{\Omega_{\delta'}}u(x,t_1)^p\,dx+t^{\frac{p}{1-m}}+ \|f\|_{L^{\infty}(\1\Omega\times (0,T))}^p\right)\quad\forall 0<t_1\le t<T.
\end{equation}
Let $C_1=\max(\|u_0\|_{L^{\infty}(\Omega)},\|u\|_{L^{\infty}(\Omega\times (0,T))})$. By the mean value theorem for any $x\in\Omega_{\delta'}$ and $0<t_1<T$ there exists a constant $\xi$ between $u(x,t_1)$ and $u_0(x)$ such that
\begin{equation*}
|u(x,t_1)^p-u_0(x)^p|=p|\xi|^{p-1}|u(x,t_1)-u_0(x)|\le pC_1^{p-1}|u(x,t_1)-u_0(x)|.
\end{equation*}
Hence
\begin{equation}\label{up-uop-integral6}
\left|\int_{\Omega_{\delta'}}u(x,t_1)^p\,dx-\int_{\Omega_{\delta'}}u_0^p\,dx\right|
\le pC_1^{p-1}\int_{\Omega_{\delta'}}|u(x,t_1)-u_0(x)|\,dx
\end{equation}
holds for any $0<t_1<T$.
Since $u$ is a solution of  \eqref{fde-Dirichlet-problem} with initial value $u_0$, letting $t_1\to 0$ in \eqref{up-uop-integral6} by Definition \ref{dirichlet-problem-soln-defn} we have,
\begin{align}\label{up-uop-integral-limit-t-0}
&\lim_{t_1\to 0}\left|\int_{\Omega_{\delta'}}u(x,t_1)^p\,dx-\int_{\Omega_{\delta'}}u_0^p\,dx\right|=0\notag\\
\Rightarrow\quad&\lim_{t_1\to 0}\int_{\Omega_{\delta'}}u(x,t_1)^p\,dx=\int_{\Omega_{\delta'}}u_0^p\,dx.
\end{align}
Since \eqref{up-fp-integral-ineqn} holds for any $t_1\in (0,t)$, letting  $t_1\to 0$ in \eqref{up-fp-integral-ineqn}, by \eqref{up-uop-integral-limit-t-0} we get \eqref{up-u0-p-integral-ineqn} and the lemma follows.

\end{proof}

By Lemma \ref{up-u0-p-integral-ineqn-lem} and an argument similar to the proof of Corollary 1.8 of \cite{H1} and a compactness argument we have the following result.

\begin{lem}\label{lem-Cor-2-2-of-Hs1}
Let $n\ge 3$, $0<m\leq\frac{n-2}{n}$, $p>\frac{n(1-m)}{2}$, $0\le f\in L^{\infty}(\partial\Omega\times (0,T))$  and $0\le u_0\in L^{\infty}(\Omega)$. Suppose $u\in L^{\infty}(\Omega\times (0,T))$ is a solution of  \eqref{fde-Dirichlet-problem}.
Then for any $0<\delta'<\delta<\delta_0$ and $0<t_1<t_2<T$ there exist  constants $C>0$ and $\theta>0$
such that
\begin{equation}\label{u-L-infty-bdary}
\left\|u\right\|_{L^{\infty}(\Omega_{\delta}\times[t_1,t_2))}\leq C\left(1+\|f\|_{L^{\infty}}^p+\int_{\Omega_{\delta'}}u_0^p\,dx\right)^{\theta/p}+\|f\|_{L^{\infty}}
\end{equation}
where $\Omega_{\delta}$, $\Omega_{\delta'}$, is given by \eqref{d-delta}.
\end{lem}

By Lemma 1.7 of \cite{H1} and  a compactness argument we have the following result.

\begin{lem}[cf. Corollary 1.8 of \cite{H1}]\label{L-infinity-initial-Lp-bd1}
Let $n\ge 3$, $0<m\leq\frac{n-2}{n}$, $p>\frac{n(1-m)}{2}$, $0\le f\in L^{\infty}(\partial\Omega\times (0,T))$  and $0\le u_0\in L^{\infty}(\Omega)$. Suppose $u\in L^{\infty}(\Omega\times (0,T))$ is a solution of  \eqref{fde-Dirichlet-problem}.
Then for any $0<\delta'<\delta<\delta_0$ and $0<t_1<t_2<T$ there exist  constants $C>0$ and $\theta>0$
such that
\begin{equation}\label{u-L-infty-interior}
\|u\|_{L^{\infty}(D_{\delta}\times [t_1,t_2))}
\le C\left(1+\int_{D_{\delta'}}u_0^p\,dx\right)^{\theta/p}
\end{equation}
where $D_{\delta}$, $D_{\delta'}$, is given by \eqref{d-delta}.
\end{lem}

\begin{thm}\label{uniqueness-thm-bd-domain}(cf. Theorem 1.1 of \cite{H2})
Let $n\geq 3$, $0<m<\frac{n-2}{n}$, $0<\delta_1<\min (1,\delta_0)$, $\mu_0>0$, $f_1, f_2\in C^3(\partial\Omega\times (0,\infty))\cap L^{\infty}(\partial\Omega\times (0,\infty))$ be such that 
\begin{equation}\label{f1-f2-ineqn}
f_2\ge f_1\ge 0\quad\mbox{ and }f_2\ge\mu_0\quad\mbox{ on }\partial\Omega\times (0,\infty)
\end{equation}
and 
\begin{equation}\label{u0-1-2-ineqn}
u_{0,1}, u_{0,2}\in L_{loc}^p(\2{\Omega}\setminus\{a_1,\cdots,a_{i_0}\}), u_{0,2}\ge u_{0,1}\ge 0, u_{0,2}\ge\mu_0\,\,\mbox{ for some constant  }p>\frac{n(1-m)}{2}.
\end{equation}
Let $\widehat{\Omega}$ be given by \eqref{omega-hat}.
Suppose  $u_1, u_2\in L_{loc}^{\infty}((\overline{\Omega}\setminus\{a_1,\dots, a_{i_0}\})\times (0,\infty))\cap C^{2,1}(\widehat{\Omega}\times (0,\infty))$ are subsolution and supersolution of \eqref{Dirichlet-blow-up-problem}  with  $f=f_1, f_2$ and $u_0=u_{0,1}, u_{0,2}$  respectively which satisfy
\begin{equation}\label{u1-2-lower-bd}
u_2(x,t)\ge\mu_0\quad\forall x\in \widehat{\Omega}, t>0
\end{equation}
such that for any constants $T>0$ and $\delta_2\in (0,\delta_1)$ there exist constants $C_1=C_1(T)>0$, $C_2=C_2(T)>0$, such that
\begin{equation}\label{u-lower-upper-blow-up-rate}
u_1(x,t)\le\frac{C_1}{|x-a_i|^{\gamma_i'}}\quad\mbox{ and }\quad u_2(x,t)\ge\frac{C_2}{|x-a_i|^{\gamma_i}} \quad \forall 0<|x-a_i|<\delta_2, 0<t<T, i=1,2,\dots,i_0.
\end{equation}
for some constants $\gamma_i$, $\gamma_i'$, $i=1,\dots,i_0$,  satisfying \eqref{gamma-gamma'-lower-bd}. Then 
\begin{equation}\label{u1-u2-compare3}
u_1(x,t)\le u_2(x,t)\quad\forall x\in \widehat{\Omega}, t>0.
\end{equation}
\end{thm}
\begin{proof}
Since the proof of the theorem is a modification of the proof of Theorem 1.1 of \cite{H2}, we will only sketch the proof here. Let 
\begin{equation*}
D_+=\{(x,t)\in \widehat{\Omega}\times (0,\infty):u_1(x,t)>u_2(x,t)\}
\end{equation*}
and $\alpha>\max(2+n,\gamma_1',\gamma_2',\dots,\gamma_{i_0}')$. Then by \eqref{u1-2-lower-bd} for any $(x,t)\in D_+$,
\begin{equation*}
u_1(x,t)>u_2(x,t)\ge\mu_0.
\end{equation*}
Hence by the mean value theorem,
\begin{equation}\label{u1-u2-difference0}
(u_1^m-u_2^m)_+(x,t)\le m\mu_0^{m-1}(u_1-u_2)_+(x,t)\quad\forall x\in \widehat{\Omega}, t>0.
\end{equation}
As in \cite{H2} we choose $\psi\in C^{\infty}(\2{\Omega}\setminus\{a_1,\cdots, a_{i_0}\})$ such that 
$\psi (x)=|x-a_i|^{\alpha}$ for any $x\in\cup_{i=1}^{i_0}B_{\delta_0}(a_i)$ and
\begin{equation}\label{psi+}
\psi (x)\ge c_1\quad\forall x\in \2{\Omega}\setminus \cup_{i=1}^{i_0}B_{\delta_2}(a_i)
\end{equation} 
for some constant $c_1>0$. Let $T>0$. Since 
$$
u_1,u_2\in L_{loc}^{\infty}((\overline{\Omega}\setminus\{a_1,\dots, a_{i_0}\})\times (0,\infty)),
$$ 
by \eqref{u-lower-upper-blow-up-rate}  and the choice of $\alpha$, we have for any $i=1,\cdots,i_0$, 
\begin{align}\label{u1-u2-blow-neighd-l1-bd}
\int_{B_{\delta_2}(a_i)}|x-a_i|^{\alpha}(u_1-u_2)_+(x,t)\,dx\le&C_T\int_0^{\delta_2}\rho^{\alpha+n-\gamma_i'-1}\,d\rho\notag\\
=&C_T'\delta_2^{\alpha+n-\gamma_i'}<\infty\quad\forall 0<t<T
\end{align}
for some constants $C_T>0$, $C_T'>0$. 
Since by the same argument as the proof of Proposition 2.2 of \cite{H2}, the result of Proposition 2.2 of \cite{H2} remains valid for $u_1$, $u_2$. That is
\begin{equation*}
\|u_i(\cdot,t)-u_{0,i}\|_{L^1(\Omega_{\delta})}\to 0\quad\mbox{ as }t\to 0\quad\forall 0<\delta<\delta_0, i=1,2.
\end{equation*} 
Hence there exists a constant $C_3(T)>0$ such that
\begin{align}\label{initial-value-l1-sense}
&\|u_i(\cdot,t)-u_{0,i}\|_{L^1(\2{\Omega}\setminus \cup_{i=1}^{i_0}B_{\delta_2}(a_i))}\le C_3(T) \quad\forall 0<t<T, i=1,2\notag\\
\Rightarrow\quad&\|u_1(\cdot,t)-u_2(\cdot,t)\|_{L^1(\2{\Omega}\setminus \cup_{i=1}^{i_0}B_{\delta_2}(a_i))}\le 2C_3(T) \quad\forall 0<t<T.
\end{align}  
By \eqref{u1-u2-blow-neighd-l1-bd} and \eqref{initial-value-l1-sense},
\begin{equation}\label{u1-u2-l1-with-psi-bd}
\int_{\widehat{\Omega}}\psi (x)(u_1-u_2)_+(x,t)\,dx\le C_T'\delta_2^{\alpha+n-\gamma_i'}+2C_3(T)<\infty\quad\forall 0<t<T.
\end{equation}
By \eqref{u-lower-upper-blow-up-rate} and the mean value theorem for any  $|x-a_i|\le\delta_2$, $0<t<T$, $i=1,\cdots,i_0$,
\begin{align}\label{u1-u2-difference-local-ineqn}
&|x-a_i|^{\alpha-2}(u_1^m-u_2^m)_+(x,t)\notag\\
\le &m|x-a_i|^{\alpha-2}u_2(x,t)^{m-1}(u_1-u_2)_+(x,t)\notag\\
\le &mC_2(T)^{m-1}|x-a_i|^{(1-m)\gamma_i-2+\alpha}(u_1-u_2)_+(x,t)\notag\\
\le &mC_2(T)^{m-1}\delta_0^{(1-m)\gamma_i-2}|x-a_i|^{\alpha}(u_1-u_2)_+(x,t)\notag\\
\le &mC_2(T)^{m-1}\delta_0^{(1-m)\gamma_i-2}\psi(x)(u_1-u_2)_+(x,t).
\end{align}
As in \cite{H2} we now choose a nonnegative monotone increasing function $\phi\in C^{\infty}(\Rr)$ such that $\phi(s)=0$ for any $s\le 1/2$ and $\phi(s)=1$ for any $s\ge 1$. For any $0<\delta<\delta_0$, let $\phi_{\delta}(x)=\phi (|x|/\delta)$ and
\begin{equation*}
w_{\delta}(x)=\Pi_{i=1}^{i_0}\phi_{\delta}(x-a_i).
\end{equation*}
Then by  \eqref{u1-u2-difference0}, \eqref{psi+}, \eqref{u1-u2-l1-with-psi-bd} and \eqref{u1-u2-difference-local-ineqn} and an argument similar to  the proof of Theorem 1.1 of \cite{H2}, 
 \begin{align}\label{u1-2-difference-in-l1-time-derivative-ineqn3}
&\frac{\1}{\1 t}\left(\int_{\widehat{\Omega}}(u_1-u_2)_+\psi w_{\delta}\,dx\right)\notag\\
\le&C\int_{\Omega\setminus\cup_{i=1}^{i_0}B_{\delta_2}(a_i)}(u_1^m-u_2^m)_+(x,t)\,dx\notag\\
&\quad +C\int_{\cup_{i=1}^{i_0}B_{\delta_2}(a_i)}|x-a_i|^{\alpha-2}(u_1^m-u_2^m)_+(x,t)\,dx\notag\\
\le&C\int_{\widehat{\Omega}}(u_1-u_2)_+(x,t)\psi(x)\,dx\notag\\
&\quad +C\int_{\cup_{i=1}^{i_0}B_{\delta_2}(a_i)}|x-a_i|^{\alpha-2}(u_1^m-u_2^m)_+(x,t)\,dx\notag\\
\le& C_T\int_{\widehat{\Omega}}(u_1-u_2)_+(x,t)\psi(x)\,dx. 
\end{align}
Integrating \eqref{u1-2-difference-in-l1-time-derivative-ineqn3} over $(0,t)$ as letting $\delta\to 0$,
\begin{equation}\label{u1-2-difference-in-l1-time-derivative-ineqn4}
\int_{\widehat{\Omega}}(u_1-u_2)_+(x,t)\psi (x)\,dx\le 
C_T\int_0^t\int_{\widehat{\Omega}}(u_1-u_2)_+(x,t)\psi(x)\,dx\,dt\quad\forall 0<t<T.
\end{equation}
By \eqref{u1-2-difference-in-l1-time-derivative-ineqn4} and the Gronwall inequality, we get \eqref{u1-u2-compare3} and the theorem follows.
\end{proof}

\begin{thm}\label{uniqueness-thm2-bd-domain}(cf. Theorem 1.2 of \cite{H2})
Let $n\geq 3$, $0<m<\frac{n-2}{n}$, $0<\delta_1<\min (1,\delta_0)$, $\mu_0>0$, $f_1, f_2\in L^{\infty}(\partial\Omega\times (0,\infty))$, $u_{0,1}, u_{0,2}\in L_{loc}^p(\2{\Omega}\setminus\{a_1,\cdots,a_{i_0}\})$ for some constant $p>\frac{n(1-m)}{2}$, be such that \eqref{f1-f2-ineqn} and \eqref{u0-1-2-ineqn} hold. 
Let $\widehat{\Omega}$ be given by \eqref{omega-hat}. Suppose $u_1, u_2\in L_{loc}^{\infty}((\overline{\Omega}\setminus\{a_1,\dots, a_{i_0}\})\times (0,\infty))\cap C^{2,1}(\widehat{\Omega}\times (0,\infty))$ are subsolution and supersolution of \eqref{Dirichlet-blow-up-problem}  with  $f=f_1, f_2$ and $u_0=u_{0,1}, u_{0,2}$ respectively which satisfy \eqref{u1-2-lower-bd} such that for any constants $T>0$ and $\delta_2\in (0,\delta_1)$ there exists a constant $C_2=C_2(T)>0$ such that 
\begin{equation}\label{u-lower-upper-blow-up-rate2}
u_j(x,t)\le\frac{C_2}{|x-a_i|^{\gamma_i'}} \quad \forall 0<|x-a_i|<\delta_2, 0<t<T, i=1,2,\dots,i_0,j=1,2
\end{equation} 
holds for some constants  
\begin{equation}\label{gamma-gamma'-lower-upper-bd2}
\frac{2}{1-m}<\gamma_i'<\frac{n-2}{m}\quad\forall i=1,2,\dots,i_0.
\end{equation}
Then \eqref{u1-u2-compare3} holds.

\end{thm}
\begin{proof}
Since the proof is similar to the proof of Theorem 1.2 of \cite{H2}, we will only sketch the argument here. Similar to the proof of Theorem 1.2 of \cite{H2} we let
\begin{equation*}
A(x,t)=\left\{\begin{aligned}
&\frac{u_1(x,t)^m-u_2(x,t)^m}{u_1(x,t)-u_2(x,t)}\quad\forall x\in\widehat{\Omega}, t>0 \mbox{ satisfying }u_1(x,t)\ne u_2(x,t)\\
&mu_2(x,t)^{m-1}\qquad\qquad\forall x\in\widehat{\Omega}, t>0 \mbox{ satisfying }u_1(x,t)= u_2(x,t)\\
&0\qquad\qquad\qquad\qquad\,\,\forall x=a_i, i=1,\cdots, i_0,t>0.
\end{aligned}\right.
\end{equation*}
For any $k\in\Z^+$, let
\begin{equation*}
\alpha_k(x,t)=\left\{\begin{aligned}
&\frac{|u_1(x,t)^m-u_2(x,t)^m|}{|u_1(x,t)-u_2(x,t)|+(1/k)}\quad\forall x\in\widehat{\Omega}, t>0\\
&0\qquad\qquad\qquad\qquad\qquad\quad\forall x=a_i, i=1,\cdots, i_0,t>0
\end{aligned}\right.
\end{equation*}
and $A_k(x,t)=\alpha_k(x,t)+k^{-1}$. We claim that the function $A(x,t)\in L^{\infty}(\widehat{\Omega}\times (0,\infty))$. We divide the  proof of  this claim into two cases.  

\noindent{\bf Case 1}: $u_2(x,t)\ge 2u_1(x,t)$.

\noindent By \eqref{u1-2-lower-bd},

\begin{equation*}
|A(x,t)|\le\frac{u_2(x,t)^m}{\frac{1}{2}u_2(x,t)}=2u_2(x,t)^{m-1}\le 2\mu_0^{m-1}.
\end{equation*}

\noindent{\bf Case 2}: $u_2(x,t)<2u_1(x,t)$.

\noindent By \eqref{u1-2-lower-bd} and the mean value theorem there exists a constant 
$\xi=\xi(x,t)$ lying between $u_1(x,t)$ and $u_2(x,t)$ such that

\begin{equation*}
|A(x,t)|\le m\xi^{m-1}\le m(u_2(x,t)/2)^{m-1}\le2^{1-m}m\mu_0^{m-1}.
\end{equation*}
By case 1 and case 2, $A(x,t)\in L^{\infty}(\widehat{\Omega}\times (0,\infty))$. Since $|\alpha_k(x,t)|\le |A(x,t)|$, we get $\alpha_k(x,t)\in L^{\infty}(\widehat{\Omega}\times (0,\infty))$ and hence  one can apply the same argument as the proof of Theorem 1.2 of \cite{H2} to conclude that the theorem holds.

\end{proof}

\begin{proof}[Proof of Theorem \ref{unique-soln-thm-bd-domain}:]Since the proof is similar to the proof of Theorem 1.1 of \cite{HK2}, we will only sketch the argument here. For any $M>0$, $0<\3<1$, let 
\begin{equation}\label{eq-def-of-u-0-epsilon-M}
\begin{cases}
\begin{aligned}
u_{0,\3}(x)&=\left(u_0(x)^m+\3^m\right)^{1/m} \\
u_{0,\3,M}(x)&=\left(\min\left(u_0(x)^m,M^m\right)+\3^m\right)^{1/m}  
\end{aligned}
\end{cases}
\end{equation}
and
\begin{equation}\label{eq-definition-of-f-epsilon-by-f}
f_{\3}(x,t)=\left(f(x,t)^m+\3^m\right)^{1/m} \qquad \forall (x,t)\in\partial\Omega\times(0,\infty).
\end{equation}
Let $u_{\3,M}$ be the solution of 
\begin{equation}\label{fde-Dirichlet-problem2}
\begin{cases}
\begin{aligned}
u_t&=\Delta u^m \qquad \mbox{in $\Omega\times(0,T)$}\\
u&=f_{\3} \qquad \quad \mbox{on $\partial\Omega\times(0,T)$}\\
u(x,0)&=u_{0,\3,M}(x) \quad\,\,\, \mbox{in $\Omega$}.
\end{aligned}
\end{cases}
\end{equation} Then
\begin{equation}\label{u-epsilon-M-mononotone}
\left\{\begin{aligned}
&u_{\3,M_2}\ge u_{\3,M_1}\ge\3\quad\,\,\mbox{ in }\Omega\times (0,\infty)\quad\forall M_2>M_1>0,\3>0\\
&u_{\3_1,M}\ge u_{\3_2,M}\qquad\quad\mbox{ in }\Omega\times (0,\infty)\quad\forall M>0,\3_1>\3_2>0.
\end{aligned}\right.
\end{equation}
By Lemma \ref{lem-Cor-2-2-of-Hs1} for any $0<\delta'<\delta<\delta_0$, $t_0'>t_0>0$, there exists a constant $C>0$ such that 
\begin{equation}\label{u-epsilon-M-local-L-infty}
\left\|u_{\3,M}\right\|_{L^{\infty}(\Omega_{\delta}\times[t_0,t_0'])}\leq C\left(1+\|f\|_{L^{\infty}}^p+\int_{\Omega_{\delta'}}u_0^p\,dx\right)^{\theta/p}+\|f\|_{L^{\infty}}=:C_0
\end{equation}
holds for any $0<\3\le 1$ and $M>0$  where $\Omega_{\delta}$, $\Omega_{\delta'}$, is given by \eqref{d-delta}.
As in \cite{HK2}, by \eqref{u-epsilon-M-mononotone} and \eqref{u-epsilon-M-local-L-infty}, as $M\to\infty$, $u_{\3,M}$ will increase monotonically to some solution $u_{\3}$ of 
\begin{equation}\label{Dirichlet-blow-up-problem2}
\left\{\begin{aligned}
u_t=&\Delta u^m\quad\,\mbox{ in }\widehat{\Omega}\times (0,\infty)\\
u=&f_{\3}\qquad\,\,\,\mbox{ on }\partial\Omega\times (0,\infty)\\
u(a_i,t)=&\infty\qquad\,\forall t>0, i=1,2\dots,i_0\\
u(x,0)=&u_{0,\3}(x)\quad\mbox{in }\widehat{\Omega}.
\end{aligned}\right.
\end{equation} 
Letting $M\to\infty$ in \eqref{u-epsilon-M-mononotone} and \eqref{u-epsilon-M-local-L-infty}, 
\begin{equation}\label{u-epsilon-positive}
\left\{\begin{aligned}
&u_{\3}\ge\3\qquad\mbox{ in }\widehat{\Omega}\times (0,\infty)\\
&u_{\3_1}\ge u_{\3_2}\quad\mbox{ in }\widehat{\Omega}\times (0,\infty)\quad\forall \3_1>\3_2>0\\
&u_{\3}\le C_0\quad\mbox{ in }\overline{\Omega_{\delta}}\times [t_0,t_0']\quad\forall 0<\3\le 1.
\end{aligned}\right.
\end{equation}
Moreover $u_{\3}$ will decrease monotonically to a solution $u$ of \eqref{Dirichlet-blow-up-problem}  as $\3\to 0$. By an argument similar to the proof of Theorem 1.1 of \cite{HK2} for any $T>0$, $\delta_2\in (0,\delta_1)$, there exists constants $C_1=C_1(T)>0$, $C_2=C_2(T)>0$, depending only on $\lambda_1$,
$\cdots$, $\lambda_{i_0}$, $\lambda_1'$, $\cdots$, $\lambda_{i_0}'$, $\gamma_1$, $\cdots$, $\gamma_{i_0}$, $\gamma_1'$, $\cdots$, $\gamma_{i_0}'$, such that both $u$ and $u_{\3}$ satisfy \eqref{u-lower-upper-blow-up-rate5} for any $0<\3<1$. 

Suppose $v$ is another solution of \eqref{Dirichlet-blow-up-problem} which satisfies \eqref{u-lower-upper-blow-up-rate5} for some constants $C_1>0$, $C_2>0$. Since by  \eqref{eq-def-of-u-0-epsilon-M} and \eqref{eq-definition-of-f-epsilon-by-f}, \begin{equation*}
u_{0,\3}\ge\max(u_0,\3)\quad\mbox{ and }f_{\3}\ge\max(f,\3),
\end{equation*} 
by  Theorem \ref{uniqueness-thm-bd-domain}, 
\begin{align*}
&v\le u_{\3}\quad\mbox{ in }\widehat{\Omega}\times (0,\infty)\quad\forall 0<\3<1\\
\Rightarrow\quad&v\le u\quad\,\mbox{ in }\widehat{\Omega}\times (0,\infty)\quad\mbox{ as }\3\to 0.
\end{align*}
Hence $u$ is the maximal solution of \eqref{Dirichlet-blow-up-problem}. 

\noindent{\bf Proof of (i) of Theorem \ref{unique-soln-thm-bd-domain}}: 

\noindent  Suppose there exist constants $T_0'>T_0>0$ and $\mu_1>0$ such that
\eqref{f-positive-after-t1} holds and $T_1\in (T_0,T_0')$.  Let $T_3=(T_0+T_1)/2$, $T=(T_1-T_0)/2$ and $C_3=T^{\frac{2}{1-m}}/\mu_1^2$. Let $q$ and $\lambda_1>0$ be the first positive eigenfunction and the first eigenvalue of $-\Delta$ on $\Omega$. By the proof of Theorem 2.2 of \cite{H1} there exists a constant $C_4>0$ such that the function
\begin{equation}\label{w-defn}
w(x,t)=\frac{[m(t-T_3)]^{1/(1-m)}}{(C_3+C_4q(x))^{1/2}}
\end{equation}
is a subsolution of  \eqref{fde} in $\Omega\times (T_3,\infty)$.  Since $w(x,T_3)=0$ in $\Omega$ and 
\begin{equation*}
w(x,t)=\left(\frac{m(t-T_3)}{T}\right)^{1/(1-m)}\mu_1\le f_{\3}(x,t)\quad\mbox{ on }\partial\Omega\times [T_3,T_1],
\end{equation*}
by Theorem \ref{uniqueness-thm-bd-domain},
\begin{align}\label{u-lower-positive-bd}
&u_{\3}(x,t)\ge w(x,t)\quad\forall x\in\widehat{\Omega}\times (T_3, T_1], 0<\3<1\notag\\
\Rightarrow\quad&u(x,t)\ge w(x,t)\quad\forall x\in\widehat{\Omega}\times (T_3, T_1]\quad\mbox{ as }\3\to 0\notag\\
\Rightarrow\quad&u(x,T_1)\ge\mu_3:=\left(\frac{m(T_1-T_0)}{2}\right)^{1/(1-m)}(C_3+C_4\|q\|_{\infty})^{-1/2}\quad\forall x\in\widehat{\Omega}
\end{align}
Let $\mu_2=\min (\mu_1,\mu_3)$.
Then by \eqref{f-positive-after-t1}, \eqref{u-lower-positive-bd} and Theorem \ref{uniqueness-thm-bd-domain},
\begin{align*}
&u_{\3}(x,t)\ge\mu_2\quad\forall x\in\widehat{\Omega}\times [T_1,T_0'), 0<\3<1\notag\\
\Rightarrow\quad&u(x,t)\ge \mu_2\quad\forall x\in\widehat{\Omega}\times [T_1,T_0')\quad\mbox{ as }\3\to 0
\end{align*}
and (i) follows.

\noindent{\bf Proof of (ii) of Theorem \ref{unique-soln-thm-bd-domain}}:
 
\noindent Suppose there exists a constant $T_2\ge 0$ such that
\eqref{eq-condition-monotone-decreasingness-of-f2-34} holds. Then $f_{\3}$ is monotone decreasing in $t$ on 
$\partial\Omega\times (T_2,\infty)$. Hence similar to Theorem 1.1 of \cite{HK2} both $u_{\3,M}$ and $u_{\3}$ 
satisfies \eqref{Aronson-Bernilan-ineqn}. Putting $u=u_{\3}$ in \eqref{Aronson-Bernilan-ineqn} and letting $\3\to 0$, we get that $u$ satisfies \eqref{Aronson-Bernilan-ineqn}
and (ii) follows.
\end{proof}

By Lemma \ref{lem-Cor-2-2-of-Hs1}, Lemma \ref{L-infinity-initial-Lp-bd1} and the construction of solution of \eqref{Dirichlet-blow-up-problem} in  Theorem \ref{unique-soln-thm-bd-domain} we recover Lemma 2.9 of \cite{HK2} and have the following results.

\begin{lem}[cf.  Lemma 2.9 of \cite{HK2}]\label{lem-Cor-2-2-of-Hs1-3}
Let $n\ge 3$, $0<m\leq\frac{n-2}{n}$, $0\le f\in L^{\infty}(\partial\Omega\times[0,\infty))$  and $0\le u_0\in L^p_{loc}(\overline{\Omega}\setminus\{a_1,\cdots,a_{i_0}\})$ for some constant $p>\frac{n(1-m)}{2}$. Suppose $u$ is a solution  of \eqref{Dirichlet-blow-up-problem}.
Then for any $0<\delta'<\delta<\delta_0$ and $0<t_1<t_2<T$ there exist  constants $C>0$ and $\theta>0$
such that \eqref{u-L-infty-bdary} holds.
\end{lem}

\begin{lem}\label{L-infinity-initial-Lp-bd1-4}
Let $n\ge 3$, $0<m\leq\frac{n-2}{n}$, $0\le f\in L^{\infty}(\partial\Omega\times[0,\infty))$  and $0\le u_0\in L^p_{loc}(\overline{\Omega}\setminus\{a_1,\cdots,a_{i_0}\})$ for some constant $p>\frac{n(1-m)}{2}$. Suppose $u$ is a solution  of \eqref{Dirichlet-blow-up-problem}.
Then for any $0<\delta'<\delta<\delta_0$ and $0<t_1<t_2<T$ there exist  constants $C>0$ and $\theta>0$
such that \eqref{u-L-infty-interior} holds.
\end{lem}

\begin{lem}\label{up-u0-p-integral-ineqn-lem-5}
Let $n\ge 3$, $0<m\leq\frac{n-2}{n}$, $0\le f\in L^{\infty}(\partial\Omega\times[0,\infty))$  and $0\le u_0\in L^p_{loc}(\overline{\Omega}\setminus\{a_1,\cdots,a_{i_0}\})$ for some constant $p>\frac{n(1-m)}{2}$. Suppose $u$ is a solution  of \eqref{Dirichlet-blow-up-problem}.
Then for any $0<\delta'<\delta<\delta_0$ there exists a  constant $C>0$ depending only on $p$, $m$, $\delta$ and $\delta'$ such that \eqref{up-u0-p-integral-ineqn} holds.
\end{lem}

\begin{rmk}
By an argument similar to the proof of Theorem \ref{unique-soln-thm-bd-domain} but with Theorem \ref{uniqueness-thm2-bd-domain} replacing Theorem \ref{uniqueness-thm-bd-domain} in the proof we get Theorem \ref{unique-soln-thm-bd-domain2}.
\end{rmk}

By Theorem \ref{uniqueness-thm-bd-domain}, Theorem \ref{uniqueness-thm2-bd-domain} and the construction of solution of \eqref{Dirichlet-blow-up-problem} in  Theorem \ref{unique-soln-thm-bd-domain} we have the following corollaries.

\begin{cor}\label{comparison-cor-bd-domain}
Let $n\geq 3$, $0<m<\frac{n-2}{n}$, $0<\delta_1<\min (1,\delta_0)$ and $T_1>0$. Let $0\le u_{0,1}\le u_{0,2}\in L_{loc}^p(\2{\Omega}\setminus\{a_1,\cdots,a_{i_0}\})$ for some constant  $p>\frac{n(1-m)}{2}$ and $f_1, f_2\in C^3(\partial\Omega\times (0,T_1))\cap L^{\infty}(\partial\Omega\times (0,T_1))$ be such that 
\begin{equation*}
f_2\ge f_1\ge 0\quad\mbox{ on }\partial\Omega\times (0,T_1)
\end{equation*}
holds. Let $\widehat{\Omega}$ be given by \eqref{omega-hat}.
Suppose  $u_1, u_2\in L_{loc}^{\infty}((\overline{\Omega}\setminus\{a_1,\dots, a_{i_0}\})\times (0,T_1))\cap C^{2,1}(\widehat{\Omega}\times (0,T_1))$ are the maximal solutions  of \eqref{Dirichlet-blow-up-problem}  in $\widehat{\Omega}\times (0,T_1)$ with $f=f_1, f_2$ and $u_0=u_{0,1}, u_{0,2}$  respectively 
such that for any constants $0<T<T_1$ and $\delta_2\in (0,\delta_1)$ there exist constants $C_1=C_1(T)>0$, $C_2=C_2(T)>0$, such that 
\begin{equation}\label{u-1-2-lower-upper-blow-up-rate5}
\frac{C_1}{|x-a_i|^{\gamma_i}}\le u_j(x,t)\le\frac{C_2}{|x-a_i|^{\gamma_i'}} \quad \forall 0<|x-a_i|<\delta_2, 0<t<T, i=1,2,\dots,i_0,j=1,2
\end{equation} holds
for some constants $\gamma_i$, $\gamma_i'$, $i=1,\dots,i_0$,  satisfying \eqref{gamma-gamma'-lower-bd}. Then \eqref{u1-u2-compare3} holds for any $x\in \widehat{\Omega}$, $0<t<T_1$. 
\end{cor}

\begin{cor}\label{comparison-cor2-bd-domain}
Let $n\geq 3$, $0<m<\frac{n-2}{n}$, $0<\delta_1<\min (1,\delta_0)$ and $T_1>0$. Let $0\le u_{0,1}\le u_{0,2}\in L_{loc}^p(\2{\Omega}\setminus\{a_1,\cdots,a_{i_0}\})$ for some constant  $p>\frac{n(1-m)}{2}$ and $0\le f_1\le f_2\in L^{\infty}(\partial\Omega\times (0,T_1))$. Let $\widehat{\Omega}$ be given by \eqref{omega-hat}.
Suppose  $u_1, u_2\in L_{loc}^{\infty}((\overline{\Omega}\setminus\{a_1,\dots, a_{i_0}\})\times (0,\infty))\cap C^{2,1}(\widehat{\Omega}\times (0,\infty))$ are the maximal solutions of \eqref{Dirichlet-blow-up-problem} in $\widehat{\Omega}\times (0,T_1)$ with  $f=f_1, f_2$ and $u_0=u_{0,1}, u_{0,2}$ respectively which satisfy \eqref{u1-2-lower-bd} such that for any constants $0<T<T_1$ and $\delta_2\in (0,\delta_1)$ 
there exist constants $C_1=C_1(T)>0$, $C_2=C_2(T)>0$, such that \eqref{u-1-2-lower-upper-blow-up-rate5} holds for some constants $\gamma_i$, $\gamma_i'$, $i=1,\dots,i_0$,  satisfying \eqref{gamma-gamma'-lower-upper-bd}. Then \eqref{u1-u2-compare3} holds for any $x\in \widehat{\Omega}$, $0<t<T_1$.  
\end{cor}

\section{Asymptotic behaviour of blow-up solutions}
\setcounter{equation}{0}
\setcounter{thm}{0}

In this section we will prove the asymptotic behaviour of the maximal finite blow-up points solutions.  

\begin{proof}[\textbf{Proof of Theorem \ref{convergence-thm-bded-domain}}:]
For any $0<\3<1$, let  $u_{0,\3}$, $f_{\3}$ and $u_{\3}$ as in the proof of Theorem \ref{unique-soln-thm-bd-domain}. Then 
\begin{equation}\label{u-u-epsilon-ineqn}
u(x,t)\le u_{\3}(x,t)\quad\forall x\in\widehat{\Omega}, t>0.
\end{equation}
By an argument similar to the proof of Theorem 1.1 of \cite{HK2} for any $T>0$, $\delta_2\in (0,\delta_1)$, there exists constants $C_1=C_1(T)>0$, $C_2=C_2(T)>0$, depending only on $\lambda_1$,
$\cdots$, $\lambda_{i_0}$, $\lambda_1'$, $\cdots$, $\lambda_{i_0}'$, $\gamma_1$, $\cdots$, $\gamma_{i_0}$, $\gamma_1'$, $\cdots$, $\gamma_{i_0}'$, such that
\eqref{u-lower-upper-blow-up-rate5} holds with $u=u_{\3}$ for all $0<\3<1$. For any $0<\delta<\delta_0$, let $\Omega_{\delta}$ be given by \eqref{d-delta}. By \eqref{gamma-gamma'-lower-upper-bd} and Lemma 3.2 of \cite{HK2} for any constants $0<\delta<\delta_0$, $t_0>0$, there exists a constant $C_{\delta}>0$ such that
\begin{equation}\label{u-epsilon-uniform-upper-bd6}
u_{\3}(x,t)\le C_{\delta}\quad\forall x\in\overline{\Omega_{\delta}}\times [t_0,\infty), 0<\3<1.
\end{equation}
Let $\{t_i\}_{i=1}^{\infty}\subset\mathbb{R}^+$ be a sequence such that $t_i\to\infty$ as $i\to\infty$. 
Let $u_i(x,t)=u(x,t+t_i)$ and $u_{\3,i}=u_{\3}(x,t+t_i)$. Let $\phi_{\3}$ be the solution of \eqref{harmonic-eqn} with $g^m$ being replaced by $g^m+\3^m$.
By  Theorem 1.5 of \cite{HK2} and \eqref{f-t-infty-limit}, \eqref{eq-def-of-u-0-epsilon-M}, \eqref{eq-definition-of-f-epsilon-by-f},
\begin{equation}\label{u-epsilon-limit}
u_{\3}\to \phi_{\3}^{\frac{1}{m}}\quad\mbox{ uniformly in }C^2(K)\quad\mbox{ as }t\to\infty
\end{equation} 
for any compact subset $K\subset \overline{\Omega}\setminus\{a_1,\dots,a_{i_0}\}$.

We now divide the proof into three cases.

\noindent{\bf Case (i)}: $g>0$ on $\partial\Omega$.

\noindent Since $\min_{\partial\Omega}g>0$, we can choose a constant $\mu_1\in (0,\min_{\partial\Omega}g)$. Then by \eqref{f-t-infty-limit} there exists a constant $T_0>0$ such that \eqref{f-positive-after-t1} holds with $T_0'=\infty$. Let $T_1>T_0$. Then by Theorem \ref{unique-soln-thm-bd-domain} there exists a constant $\mu_2\in (0,\mu_1)$ such that
\eqref{u-eventual-positive} holds with $T_0'=\infty$. By 
Theorem 1.5 of \cite{HK2}, \eqref{f-t-infty-limit}, \eqref{f-positive-after-t1} and  \eqref{u-eventual-positive}, \eqref{u-phi-1/m-limit1} holds
for any compact subset $K$ of $\overline{\Omega}\setminus\{a_1,\dots,a_{i_0}\}$ and (i) follows.

\noindent{\bf Case (ii)}: $g\not\equiv 0$ on $\partial\Omega$ and \eqref{f>g}, \eqref{u0-lower-bd3} holds.

\noindent Since $f_{\3}\ge\max(g,\3)$ and the function $\phi^{1/m}\in C^1(\overline{\Omega})\cap C^{2,1}(\Omega)$ satisfy 
\eqref{fde-Dirichlet-problem} with $f=g$ and $u_0=\phi$, by \eqref{f>g}, \eqref{u0-lower-bd3}  and Theorem \ref{uniqueness-thm2-bd-domain},
\begin{align}\label{u-epsilon-lower-bd6}
&u_{\3}(x,t)\ge \phi(x)^{1/m}\quad\forall x\in \widehat{\Omega}, t>0, 0<\3<1\notag\\
\Rightarrow\quad &u(x,t)\ge \phi(x)^{1/m}\quad\forall x\in \widehat{\Omega}, t>0\quad\mbox{ as } \3\to 0.
\end{align}
Since $\phi(x)>0$ on $\Omega$, by  \eqref{u-u-epsilon-ineqn}, \eqref{u-epsilon-uniform-upper-bd6} and \eqref{u-epsilon-lower-bd6} for any $N>0$ the equation \eqref{fde} for the sequence $\{u_i\}_{t_i>-N}$ is uniformly parabolic on any compact subset of $K\subset \widehat{\Omega}\times [-N,N]$. Hence by the parabolic Schauder estimates \cite{LSU} the sequence $\{u_i\}_{t_i>-N}$ is uniformly continuous in $C^2(K)$ for any compact subset of $K\subset \widehat{\Omega}\times [-N,N]$. Thus by \eqref{u-u-epsilon-ineqn}, \eqref{u-epsilon-limit}, \eqref{u-epsilon-lower-bd6},  the Ascoli Theorem and a diagonalization argument the sequence $\{u_i\}$ has a subsequence which we may assume without loss of generality to be the sequence itself that converges uniformly in $C^2(K)$ for any compact subset of $K\subset \widehat{\Omega}\times (-\infty,\infty)$ to a  solution $v$ of \eqref{fde} in 
$\widehat{\Omega}\times(-\infty,\infty)$ which satisfies
\begin{align*}
&\phi(x)^{1/m}\le v(x,t)\le \phi_{\3}(x)^{1/m}\quad\forall x\in \widehat{\Omega}, t>0\notag\\
\Rightarrow\quad&v(x,t)=\phi(x)^{1/m}\qquad\qquad\quad\,\,\,\forall x\in \widehat{\Omega}, t>0\quad\mbox{ as }\3\to 0\notag\\
\Rightarrow\quad&u(x,t_i)\to v(x,0)=\phi(x)^{1/m}\quad\mbox{uniformly on }C^2(K)\quad\mbox{ as }i\to\infty
\end{align*}
for any compact subset $K\subset \widehat{\Omega}$.
Since the sequence $\{t_i\}$ is arbitrary, we get \eqref{u-phi-1/m-limit1} and (ii) follows.

\noindent{\bf Case (iii)}: $g=0$ on $\partial\Omega$.

\noindent By  \eqref{u-u-epsilon-ineqn}, \eqref{u-epsilon-uniform-upper-bd6} and Theorem 1.1 of \cite{S} for any $N>0$  the sequence $\{u_i\}_{t_i>-N}$ is uniformly  continuous in $K$ for any compact subset of $K\subset \widehat{\Omega}\times [-N,N]$. Thus by \eqref{u-u-epsilon-ineqn}, \eqref{u-epsilon-uniform-upper-bd6} , Theorem 1.1 of \cite{S}, the Ascoli Theorem and a diagonalization argument the sequence $\{u_i\}$ has a subsequence which we may assume without loss of generality to be the sequence itself that converges uniformly in $K$ for any compact subset of $K\subset \widehat{\Omega}\times (-\infty,\infty)$ to a  continuous function $v$ which satisfies
\begin{align*}
&0\le v(x,t)\le\psi_{\3}(x)^{1/m}\quad\forall x\in\widehat{\Omega}, -\infty<t<\infty\notag\\
\Rightarrow\quad &v(x,t)=0\qquad\qquad\quad\quad\forall x\in\widehat{\Omega}\quad\mbox{ as }\3\to 0.
\end{align*}
Hence 
\begin{equation*}
u(x,t_i)=u_i(x,0)\to 0\quad\mbox{ as }i\to\infty.
\end{equation*}
Since the sequence $\{t_i\}$ is arbitrary, we get \eqref{u-phi-1/m-limit2} and (iii) follows.

\end{proof}

By an argument similar to the proof of Theorem \ref{convergence-thm-bded-domain} but with Remark 3.7 of \cite{HK2} replacing Theorem 1.5 of \cite{HK2} in the argument Theorem \ref{convergence-thm2-bded-domain} follows.

\section{Existence of finite blow-up solutions that blow-up at the lateral boundary}
\setcounter{equation}{0}
\setcounter{thm}{0}

In this section we will construct a solution of \eqref{Dirichlet-blow-up-problem} with appropriate  lateral boundary value such that the finite blow-up points  solution will oscillate between two given harmonic functions as $t\to\infty$. We will also prove the existence  of finite blow-up solutions that blow-up at the lateral boundary of the bounded cylindrical domain.

\begin{proof}[\textbf{Proof of Theorem \ref{oscillation-thm-bded-domain}}:]
Let $f_1=g_1$ and $u_1$ be the maximal solution of \eqref{Dirichlet-blow-up-problem} given by Theorem \ref{unique-soln-thm-bd-domain2} with $f=f_1$. 
For any $0<\delta<\delta_0$, let $D_{\delta}$ be given by \eqref{d-delta}. Let $t_0=0$ and $\delta_k=\delta_1/k$ for any $k\in\Z^+$. Then by Theorem \ref{convergence-thm2-bded-domain} there exists a constant $t_1>0$ such that 
\begin{equation}
|u_1(x,t)-\phi_1(x)|<1\quad\forall x\in D_{\delta_1},t\ge t_1.
\end{equation}
Let $f_2(x,t)=g_1(x)$ for $0<t\le t_1$ and $f_2(x,t)=g_2(x)$ for $t>t_1$. Let  
$u_2$ be the maximal solution of \eqref{Dirichlet-blow-up-problem} with $f=f_2$.
Then by Theorem \ref{convergence-thm2-bded-domain}, there exists a constant $t_2>t_1+1$ such that
\begin{equation}
|u_2(x,t)-\phi_2(x)|<\frac{1}{2}\quad\forall x\in D_{\delta_2},t\ge t_2.
\end{equation}
By repeating the above argument there exist sequences $\{t_i\}_{i=1}^{\infty}$, $t_i+1<t_{i+1}$ for all 
$i\in\Z^+$, $\{f_i\}_{i=1}^{\infty}\subset L^{\infty}(\partial\Omega)$, such that $\forall i\in\Z^+$,
\begin{equation}\label{fi-odd-defn}
f_{2i+1}(x,t)=\left\{\begin{aligned}
&g_1(x)\quad\forall x\in\partial\Omega, t\in\cup_{k=0}^{i-1}(t_{2k},t_{2k+1}]\cup (t_{2i},\infty)\\
&g_2(x)\quad\forall x\in\partial\Omega, t\in\cup_{k=1}^i(t_{2k-1},t_{2k}]
\end{aligned}\right.
\end{equation}
and 
\begin{equation}\label{fi-even-defn}
f_{2i}(x,t)=\left\{\begin{aligned}
&g_2(x)\quad\forall x\in\partial\Omega, t\in\cup_{k=1}^{i-1}(t_{2k-1},t_{2k}]\cup (t_{2i-1},\infty)\\
&g_1(x)\quad\forall x\in\partial\Omega, t\in\cup_{k=1}^i(t_{2k-2},t_{2k-1}]
\end{aligned}\right.
\end{equation}
and a sequence $\{u_i\}_{i=1}^{\infty}$ of maximal solutions of \eqref{Dirichlet-blow-up-problem} with $f=f_i$ 
that satisfies
\begin{equation}\label{ui-asymptotic-beaviour}
\left\{\begin{aligned}
&|u_{2i+1}(x,t)-\phi_1(x)|<\frac{1}{2i+1}\quad\forall x\in D_{\delta_{2i+1}},t\ge t_{2i+1}, i\in\Z^+\\
&|u_{2i}(x,t)-\phi_2(x)|<\frac{1}{2i}\qquad\quad\,\,\forall x\in D_{\delta_{2i}},t\ge t_{2i}, i\in\Z^+.
\end{aligned}\right.
\end{equation}
Let $u$ be the maximal solution of \eqref{Dirichlet-blow-up-problem} with 
\begin{equation}\label{f-defn}
f(x,t)=\left\{\begin{aligned}
&g_1(x)\quad\forall x\in\partial\Omega, t\in\cup_{k=0}^{\infty}(t_{2k},t_{2k+1}]\\
&g_2(x)\quad\forall x\in\partial\Omega, t\in\cup_{k=1}^{\infty}(t_{2k-1},t_{2k}].
\end{aligned}\right.
\end{equation}
Then by \eqref{fi-odd-defn}, \eqref{fi-even-defn} and \eqref{f-defn},
\begin{equation}\label{f=fi}
f(x,t)=f_i(x,t)\quad\forall x\in\partial\Omega, t\in (0,t_i), i\in\Z^+.
\end{equation}
Hence by Corollary \ref{comparison-cor2-bd-domain},
\begin{equation}\label{u=ui}
u(x,t)=u_i(x,t)\quad\forall x\in\partial\Omega, t\in (0,t_i], i\in\Z^+.
\end{equation}
By \eqref{ui-asymptotic-beaviour} and \eqref{u=ui},
\begin{equation*}
\left\{\begin{aligned}
&|u(x,t_{2i+1})-\phi_1(x)|<\frac{1}{2i+1}\quad\forall x\in D_{\delta_{2i+1}}, i\in\Z^+\\
&|u(x,t_{2i})-\phi_2(x)|<\frac{1}{2i}\qquad\quad\,\,\forall x\in D_{\delta_{2i}}, i\in\Z^+.
\end{aligned}\right.
\end{equation*}
Since $D_{\delta_i}\subset D_{\delta_{i+1}}$ for all $i\in\Z^+$ and $\widehat{\Omega}=\cup_{i=1}^{\infty}D_{\delta_i}$, for any $0<\3<1$ and compact subset $K$ of $\widehat{\Omega}$ there exists
$k_0\in\Z^+$, $k_0>\3^{-1}$, such that 
\begin{equation*}
K\subset D_{\delta_{k_0}}\subset D_{\delta_i} \quad\forall i\ge k_0.
\end{equation*}
Hence 
\begin{equation*}
\left\{\begin{aligned}
&|u(x,t_{2i+1})-\phi_1(x)|<\3 \quad\forall x\in K, i\ge k_0\\
&|u(x,t_{2i})-\phi_2(x)|<\3\quad\,\,\,\forall x\in K, i\ge k_0
\end{aligned}\right.
\end{equation*}
and the theorem follows.

\end{proof}

\begin{proof}[\textbf{Proof of Theorem \ref{unique-soln-thm3-bd-domain}}:]

For any $0<\delta<\delta_0$, let $D_{\delta}$ be given by \eqref{d-delta}.
For any $k\in\Z^+$, let $u_k$ be the maximal solution of \eqref{Dirichlet-blow-up-problem} with $f=k$ given by Theorem \ref{unique-soln-thm-bd-domain} which satisfies \eqref{Aronson-Bernilan-ineqn} with $T_2=0$. 
By Lemma \ref{L-infinity-initial-Lp-bd1-4} and Corollary \ref{comparison-cor-bd-domain}, for any $0<\delta<\min (1,\delta_0)$, $t_0'>t_0>0$, there exists a constant $C_{\delta}>0$ such that
\begin{equation}\label{uk-upper-lower-bd}
\left\{\begin{aligned}
&u_1(x,t)\le u_k(x,t)\le u_{k+1}(x,t)\quad\forall x\in (\overline{\Omega}\setminus\{a_1,\dots,a_{i_0}\})\times (0,\infty),k\in\Z^+\\
&u_k(x,t)\le C_{\delta}\qquad\qquad\qquad\quad\,\,\,\forall x\in D_{\delta},t_0\le t\le t_0',k\in\Z^+.
\end{aligned}\right.
\end{equation}
By Theorem \ref{unique-soln-thm-bd-domain}, for any $T>0$ and $\delta_2\in (0,\delta_1)$ there exist a constant $C_1=C_1(T)>0$ such that
\begin{equation}\label{u1-lower-blow-up-bd}
u_1(x,t)\ge\frac{C_1}{|x-a_i|^{\gamma_i}}\quad \forall 0<|x-a_i|<\delta_2, 0<t<T, i=1,2,\dots,i_0.
\end{equation}
On the other hand by the proof of Lemma 2.3 of \cite{HK2} there exists  a constant $A_0>0$ such that
\begin{equation}\label{uk-blow-up-bd}
u_k(x,t)\le \frac{A_0(1+t)^{\frac{1}{1-m}}}{|x-a_i|^{\gamma_i'}(\delta_1-|x-a_i|)^{\frac{2}{1-m}}}\quad\forall 0<|x-a_i|<\delta_1,i=1,\dots,i_0, k\in\Z^+.
\end{equation}
Let $\widehat{\Omega}$ be given by \eqref{omega-hat}. 
By \eqref{uk-upper-lower-bd} the equation \eqref{fde} for the sequence $\{u_k\}_{k=1}^{\infty}$ is uniformly parabolic on any compact subset $K$ of $\widehat{\Omega}\times (0,\infty)$. Hence by the parabolic Schauder estimates \cite{LSU} the sequence $\{u_k\}_{k=1}^{\infty}$ is uniformly continuous in $C^2(K)$ for  any compact subset $K\subset \widehat{\Omega}\times (0,\infty)$. Thus by \eqref{uk-upper-lower-bd}, \eqref{u1-lower-blow-up-bd}, \eqref{uk-blow-up-bd},  the Ascoli Theorem and a diagonalization argument the sequence $\{u_k\}_{k=1}^{\infty}$ has a subsequence which we may assume without loss of generality to be the sequence itself that increases and converges uniformly in $C^2(K)$ for any compact subset  $K\subset  \widehat{\Omega}\times (0,\infty)$ to a  solution $u$ of \eqref{fde} in 
$\widehat{\Omega}\times(0,\infty)$ which satisfies
\begin{equation}\label{u-lower-blow-up-bd}
u(x,t)\ge\frac{C_1}{|x-a_i|^{\gamma_i}}\quad \forall 0<|x-a_i|<\delta_2, 0<t<T, i=1,2,\dots,i_0,
\end{equation}
\begin{equation}\label{u-blow-up-bd}
u(x,t)\le \frac{A_0(1+t)^{\frac{1}{1-m}}}{|x-a_i|^{\gamma_i'}(\delta_1-|x-a_i|)^{\frac{2}{1-m}}}\quad\forall 0<|x-a_i|<\delta_1,i=1,\dots,i_0
\end{equation}
and
\begin{equation}\label{uk-less-than-u}
u_k\le u\quad\mbox{ in }\widehat{\Omega}\times (0,\infty)\quad\forall k\in\Z^+
\end{equation}
and $u$ also satisfies \eqref{Aronson-Bernilan-ineqn} with $T_2=0$. 

By \eqref{u-lower-blow-up-bd} and \eqref{u-blow-up-bd} $u$ satisfies \eqref{u-lower-upper-blow-up-rate5} for some constants $C_1>0$, $C_2>0$.
Now by (i) of Theorem \ref{unique-soln-thm-bd-domain}, Corollary \ref{comparison-cor-bd-domain} and \eqref{uk-upper-lower-bd}, for any $T_0>0$ there exists a constant $0<\mu_{T_0}<1$ such that
\begin{equation}\label{uk-lower-bd7}
u_k(x,t)\ge u_1(x,t)\ge\mu_{T_0}\quad\forall x\in (\overline{\Omega}\setminus\{a_1,\dots,a_{i_0}\})\times [T_0,\infty),k\in\Z^+. 
\end{equation} 
By Lemma \ref{lem-Cor-2-2-of-Hs1-3} and \eqref{uk-lower-bd7}, for any $k\in\Z^+$ the equation \eqref{fde} for $u_k$ is uniformly parabolic on any compact subset  $K\subset  (\overline{\Omega}\setminus\{a_1,\dots,a_{i_0}\})\times (0,\infty)$. Hence by the parabolic Schauder estimates \cite{LSU}, $u_k\in C^{2,1}((\overline{\Omega}\setminus\{a_1,\dots,a_{i_0}\})\times (0,\infty))$ for any $k\in\Z^+$. Thus
\begin{align}\label{u-bdary=infty}
&\underset{\substack{(y,s)\to (x,t)\\
(y,s)\in\Omega\times (0,\infty)}}{\lim}u(y,s)\ge\underset{\substack{(y,s)\to (x,t)\\
(y,s)\in\Omega\times (0,\infty)}}{\lim}u_k(y,s)=k
\quad\forall (x,t)\in\partial\Omega\times (0,\infty)\notag\\
\Rightarrow\quad&\underset{\substack{(y,s)\to (x,t)\\
(y,s)\in\Omega\times (0,\infty)}}{\lim}u(y,s)=\infty
\qquad\qquad\qquad\qquad\,\,\forall (x,t)\in\partial\Omega\times (0,\infty)\quad\mbox{ as }k\to\infty.
\end{align}
We will now show that $u$ has initial value $u_0$. Since $u_k$ has initial value $u_0$ for all $k\in\Z^+$, by Lemma 3.1 of \cite{HP} and a compactness argument for any $0<\delta<\delta_0$ there exists a constant $C>0$ depending on $\delta$ such that
\begin{align}\label{u-u1-difference-L1}
&\int_{D_{\delta}}(u_k(x,t)-u_1(x,t))\,dx\le Ct^{1/(1-m)}\quad\forall t>0,k\in\Z^+\notag\\
\Rightarrow\quad&\int_{D_{\delta}}(u(x,t)-u_1(x,t))\,dx\le Ct^{1/(1-m)}\quad\forall t>0\quad\mbox{ as }k\to\infty. 
\end{align}
Hence by \eqref{u-u1-difference-L1},
\begin{align}\label{u-u0-t-to-0}
\int_{D_{\delta}}|u(x,t)-u_0(x)|\,dx\le&\int_{D_{\delta}}(u(x,t)-u_1(x,t))\,dx+\int_{D_{\delta}}|u_1(x,t)-u_0(x)|\,dx\notag\\
\le&Ct^{1/(1-m)}+\int_{D_{\delta}}|u_1(x,t)-u_0(x)|\,dx\quad\forall t>0\notag\\
\Rightarrow\quad\lim_{t\to 0}\int_{D_{\delta}}|u(x,t)-u_0(x)|\,dx=&0\quad\forall 0<\delta<\delta_0.
\end{align}
Hence by \eqref{u-lower-upper-blow-up-rate5}, \eqref{u-bdary=infty} and \eqref{u-u0-t-to-0}, $u$ is a solution of \eqref{Dirichlet-blow-up-origin-bdary-problem}. 

We will now use a modification of the proof of Lemma 2.9 of \cite{H1} and Theorem 1.1 of \cite{H2} to show that $u$ is the minimal solution of \eqref{Dirichlet-blow-up-origin-bdary-problem}.
Suppose $v$ is another solution of \eqref{Dirichlet-blow-up-origin-bdary-problem}.
Let $0<\psi\in C^{\infty}(\2{\Omega}\setminus\{a_1,\cdots, a_{i_0}\})$ be as in the proof of Theorem \ref{uniqueness-thm-bd-domain} and let $E_j$ be given by \eqref{ej}.
 Let $T>t_1>0$ and $k\in\Z^+$.  By Lemma \ref{lem-Cor-2-2-of-Hs1-3}  $u_k\in L_{loc}^{\infty}((\2{\Omega}\setminus\{a_1,\dots,a_{i_0}\}\times (0,\infty))$.  Hence by a compactness argument there exists $j_0\in\Z^+$, $j_0>1/\delta_0$, such that
\begin{equation}\label{v-uk-bdary-ineqn}
\underset{\substack{x\in\Omega\setminus E_j\\t_1\le t\le T}}{\inf}v(x,t)>\underset{\substack{x\in\Omega\setminus E_j\\t_1\le t\le T}}{\sup}u_k(x,t)\quad\forall j\ge j_0
\end{equation}
where $E_j$ is given by \eqref{ej}.
Since $v>0$ in $\widehat{\Omega}\times (0,\infty)$ and satisfies \eqref{u-lower-upper-blow-up-rate5},
\begin{equation}\label{v-lower-bd}
v(x,t)\ge\overline{\mu}_{j_0}\quad\forall x\in E_{j_0},t_1\le t\le T
\end{equation}
for some constant $\overline{\mu}_{j_0}>0$. Hence by \eqref{uk-lower-bd7},  \eqref{v-uk-bdary-ineqn} and \eqref{v-lower-bd},  
\begin{equation}\label{v-lower-bd2}
v(x,t)\ge\min(\overline{\mu}_{j_0},\mu_{t_1})>0\quad\forall x\in\widehat{\Omega},t_1\le t\le T.
\end{equation}
By \eqref{v-lower-bd2} and an argument similar to the proof of Theorem 1.1 of \cite{H2} and the proof of Theorem \ref{uniqueness-thm-bd-domain} we get
\begin{equation}\label{uk-v-integral-ineqn2}
\int_{E_j}(u_k-v)_+(x,t)\psi(x)\,dx
\le\int_{E_j}(u_k-v)_+(x,t_1)\psi(x)\,dx+C\int_{t_1}^t\int_{E_j}(u_k-v)_+\psi (x)\,dx\,dt
\end{equation}
where $C>0$ is some constant for any $t_1\le t<T$, $j\ge j_0$. By \eqref{uk-v-integral-ineqn2}
and the Gronwall inequality,
\begin{equation}\label{uk-v-integral-ineqn3}
\int_{E_j}(u_k-v)_+(x,t)\psi(x)\,dx
\le \frac{e^{Ct}}{C}\int_{E_j}(u_k-v)_+(x,t_1)\psi(x)\,dx\quad\forall t_1\le t<T, j\ge j_0.
\end{equation}
Letting $j\to\infty$ in \eqref{uk-v-integral-ineqn2},
\begin{equation}\label{uk-v-integral-ineqn4}
\int_{\widehat{\Omega}}(u_k-v)_+(x,t)\psi(x)\,dx
\le \frac{e^{Ct}}{C}\int_{\widehat{\Omega}}(u_k-v)_+(x,t_1)\psi(x)\,dx\quad\forall t_1\le t<T.
\end{equation}
We now fix $0<\delta_2<\delta_1$. Then by \eqref{u-lower-upper-blow-up-rate5} and \eqref{v-uk-bdary-ineqn}, $\forall 0<\delta<\delta_2, j>j_0$,
\begin{align}\label{uk-v-difference-integral-ineqn10}
&\int_{\widehat{\Omega}}(u_k-v)_+(x,t_1)\psi(x)\,dx\notag\\
\le&\int_{\cup_{i=1}^{i_0}B_{\delta}(a_i)}u_k(x,t_1)\psi(x)\,dx+\int_{E_j\setminus \cup_{i=1}^{i_0}B_{\delta}(a_i)}|u_k(x,t_1)-u_0(x)|\psi(x)\,dx\notag\\
&\qquad +\int_{E_j\setminus \cup_{i=1}^{i_0}B_{\delta}(a_i)}|v(x,t_1)-u_0(x)|\psi(x)\,dx
\quad\forall 0<\delta<\delta_1\notag\\
\le&C_2\sum_{i=1}^{i_0}\int_{\cup_{i=1}^{i_0}B_{\delta}(a_i)}|x-a_i|^{\alpha -\gamma_i'}\,dx+\int_{E_j\setminus \cup_{i=1}^{i_0}B_{\delta}(a_i)}|u_k(x,t_1)-u_0(x)|\psi(x)\,dx\notag\\
&\qquad +\int_{E_j\setminus \cup_{i=1}^{i_0}B_{\delta}(a_i)}|v(x,t_1)-u_0(x)|\psi(x)\,dx\notag\\
\le&C'\delta^{\alpha+n-\gamma_i'}+\int_{E_j\setminus \cup_{i=1}^{i_0}B_{\delta}(a_i)}|u_k(x,t_1)-u_0(x)|\psi(x)\,dx+\int_{E_j\setminus \cup_{i=1}^{i_0}B_{\delta}(a_i)}|v(x,t_1)-u_0(x)|\psi(x)\,dx
\end{align}
for some  constant $C'>0$. Letting first $t_1\to 0$  and then $\delta\to 0$, $j\to\infty$, in \eqref{uk-v-difference-integral-ineqn10},
\begin{equation}\label{uk-v-difference-integral-limit=0}
\lim_{t_1\to 0}\int_{\widehat{\Omega}}(u_k-v)_+(x,t_1)\psi(x)\,dx=0.
\end{equation}
Letting $t_1\to 0$ in \eqref{uk-v-integral-ineqn4}, by \eqref{uk-v-difference-integral-limit=0},
\begin{align*}
&\int_{\widehat{\Omega}}(u_k-v)_+(x,t)\psi(x)\,dx=0\quad\forall 0<t<T, k\in\Z^+\notag\\
\Rightarrow\quad&u_k(x,t)\le v(x,t)\quad\forall x\in\widehat{\Omega}, 0<t<T, k\in\Z^+\\
\Rightarrow\quad&u(x,t)\le v(x,t)\quad\forall x\in\widehat{\Omega}, 0<t<T\quad\mbox{ as }k\to\infty.
\end{align*}
Since $T>0$ is arbitrary,
\begin{equation*}
u(x,t)\le v(x,t)\quad\forall x\in\widehat{\Omega}, t>0.
\end{equation*}
Hence $u$ is the minimal solution of \eqref{Dirichlet-blow-up-origin-bdary-problem}.

We will now prove that \eqref{u-t-inftyu-infty} holds. We choose $\widetilde{\gamma}_1, \dots, \widetilde{\gamma}_{i_0}$, such that
\begin{equation*}
\frac{2}{1-m}<\widetilde{\gamma}_i<\min\left(\frac{n-2}{m},\gamma_i\right)\quad\forall i=1,\dots,i_0.
\end{equation*}
and let
\begin{equation*}
\widetilde{u}_0(x)=\left\{\begin{aligned}
&u_0(x)\qquad\qquad\quad\forall x\in\Omega\setminus\cup_{i=1}^{i_0}B_{\delta_1}(a_i)\\
&\lambda_i|x-a_i|^{-\widetilde{\gamma}_i}\quad\forall 0<|x-a_i|<\delta_1,i=1,\dots,i_0.
\end{aligned}\right.
\end{equation*}
Then
\begin{equation*}
\widetilde{u}_0(x)\le u_0(x)\quad\mbox{ in }\widehat{\Omega}.
\end{equation*}
For any $k\in\Z^+$, let $v_k$ be the maximal solution of
\begin{equation*}
\left\{\begin{aligned}
u_t=&\Delta u^m\quad\,\mbox{ in }\widehat{\Omega}\times (0,\infty)\\
u=&k\qquad\,\,\,\mbox{ on }\partial\Omega\times (0,\infty)\\
u(a_i,t)=&\infty\qquad\,\forall t>0, i=1,2\dots,i_0\\
u(x,0)=&\widetilde{u}_0(x)\qquad\mbox{ in }\widehat{\Omega}
\end{aligned}\right.
\end{equation*}
given by Theorem \ref{unique-soln-thm-bd-domain}. Then by Corollary \ref{comparison-cor-bd-domain},
\begin{equation}\label{vk-uk-compare}
v_k\le u_k\quad\mbox{ in }\widehat{\Omega}\times (0,\infty)\quad\forall k\in\Z^+.
\end{equation}
Since by Theorem \ref{convergence-thm-bded-domain},
\begin{equation}
v_k(x,t)\to k\quad\mbox{ uniformly on }\Omega_{\delta}\quad\mbox{ as }t\to\infty
\end{equation}
for any $0<\delta<\delta_0$ and $k\in\Z^+$, by \eqref{uk-less-than-u} and \eqref{vk-uk-compare},
\begin{equation}\label{u-t-infty-limit}
\underset{\substack{x\in\Omega_{\delta}\\t\to\infty}}{\liminf}\,u(x,t)\ge k\quad\forall k\in\Z^+.
\end{equation}
for any  $0<\delta<\delta_0$. Letting $k\to\infty$ in \eqref{u-t-infty-limit}, we get \eqref{u-t-inftyu-infty} and the theorem follows.

\end{proof}

\end{document}